\documentclass[12pt,a4paper,reqno]{amsart}
\usepackage{geometry} 
\usepackage{amsmath}
\usepackage{amsthm}
\usepackage{amsfonts,amssymb}
\usepackage{fancyhdr} 
\usepackage{mathrsfs}  
\usepackage{graphicx}
\usepackage[all]{xy}

\usepackage{times}
\usepackage{cases}
\usepackage{color}

\usepackage[colorlinks,linkcolor=blue,anchorcolor=blue,citecolor=blue]{hyperref}

\usepackage{float} 
\usepackage{subfig}
\usepackage{amsthm}

 \newtheorem{thm}{Theorem}[section]
 \newtheorem{lemma}[thm]{Lemma}
 \newtheorem{prop}[thm]{Proposition}
 \newtheorem{cor}[thm]{Corollary}

 \newtheorem{rmk}{Remark}[section]
 
 \theoremstyle{definition}
 \newtheorem{definition}[thm]{Definition}

\newcommand{\IH}{{\mathbb H}}

\newcommand{\IR}{{\mathbb R}}

\newcommand{\CA}{{\mathcal A}}
\newcommand{\CB}{{\mathcal B}}

\newcommand{\CL}{{\mathcal L}}

\newcommand{\CT}{{\mathcal T}}

\newcommand \di{\,\mathrm{d}}

\newcommand{\red}{\textcolor{red}}
\allowdisplaybreaks[4]
\begin{document}

\title{Rigidity, volume and angle structures of $1$-$3$ type hyperbolic polyhedral 3-manifolds}
\author{Ke Feng, Huabin Ge, Chunlei Liu}
\address{}
\curraddr{}
\email{}
\thanks{}
\date{\today}
\dedicatory{}

\begin{abstract}
In this paper, we study the rigidity of hyperbolic polyhedral 3-manifolds and the volume optimization program of angle structures. We first study the rigidity of decorated $1$-$3$ type hyperbolic polyhedral metrics on 3-manifolds which are isometric gluing of decorated $1$-$3$ type hyperbolic tetrahedra. Here a $1$-$3$ type hyperbolic tetrahedron is a truncated hyperbolic tetrahedron with one hyperideal vertex and three ideal vertices. A decorated $1$-$3$ type polyhedron is a $1$-$3$ type hyperbolic polyhedron with a horosphere centered at each ideal vertex. We show that a decorated $1$-$3$ type hyperbolic polyhedral metric is determined up to isometry and change of decorations by its curvature. We also prove several results on the volume optimization program of Casson and Rivin, i,e. Casson-Rivin's volume optimization program is shown to be still valid for $1$-$3$ type ideal triangulated 3-manifolds. We also get a strongly $1$-efficiency triangulation when assuming the existence of an angle structure. On the whole, we follow the spirit of Luo-Yang \cite{Luo2018} to prove our main results. The main differences come from that the hyperbolic tetrahedra considered here have completely different geometry with those considered in \cite{Luo2018}.
\end{abstract}

\maketitle

\tableofcontents

\section{Introduction}
\label{section-introduction}
Hyperbolic tetrahedra are the fundamental building blocks to construct hyperbolic 3-manifolds. Indeed, with the help of Thurston's gluing equation \cite{Thurston2022}, the SnapPea written by Weeks \cite{Weeks1980, CDGW2016} in the 1980s shows how to construct a cusped hyperbolic 3-manifold by ideal tetrahedra. On the contrary, the question that whether every hyperbolic 3-manifold can be geometrically decomposed into hyperbolic tetrahedra is still open, we refer \cite{Feng2022} for some recent progress. Due to the connection between hyperbolic polyhedrons and the structure of hyperbolic 3-manifolds, we are interested in the geometry of 3 dimensional spaces which are isometric gluing of hyperbolic tetrahedra in hyperbolic spaces. Particularly, we will study the rigidity and volume maximization of these spaces in this paper and the following \cite{Feng2023} and \cite{FGM}.

In Klein's projective model for the hyperbolic space, ${\IH}^3$ is identified to the open unit ball in $\mathbb{R}^3\subset \mathbb{R}P^3$.
Let $P$ be the intersection of ${\IH}^3$ with a projective tetrahedron with one vertex outside of $\overline{{\IH}^3}$, three vertices on $\partial{\IH}^3$, and all six edges intersect with ${\IH}^3$. We call a vertex on $\partial{\IH}^3$ an ideal vertex, and the vertex outside $\overline{{\IH}^3}$ a hyperideal vertex (it is called a strictly hyperideal vertex in \cite{Bao2002, Schlenker2002}). At the hyperideal vertex, there is a unique way (see \cite{Bao2002, Schlenker2002}) to truncate $P$ so that the truncated face is perpendicular to the three edges connecting the hyperideal vertex and the three ideal vertices. The truncated tetrahedron, still denoted by $P$, is called a \emph{$1$-$3$ type hyperbolic tetrahedron}. Moreover, a \emph{decorated $1$-$3$ type hyperbolic tetrahedron} is a $1$-$3$ type hyperbolic tetrahedron with a horosphere, which is called a \emph{decoration}, centered at each ideal vertex. Similarly, we define a decorated $r$-$s$ type ($r+s=4$) hyperbolic tetrahedron which has $r$ hyperideal vertices with truncations and $s$ ideal vertices with decorations.

Our first result concern the rigidity of hyperbolic polyhedral 3-manifolds produced by isometric gluing of decorated $1$-$3$ type hyperbolic tetrahedra. Let $(M,\CT)$ be a triangulated compact pseudo 3-manifold, and assume that each tetrahedron in the triangulation $\CT$ is combinatorially equivalent to a $1$-$3$ type hyperbolic tetrahedron. Such $\CT$ is called a \emph{$1$-$3$ type triangulation}. A \emph{decorated $1$-$3$ type hyperbolic polyhedral metric} on $(M,\CT)$ is obtained by replacing each tetrahedron in $\CT$ by a decorated $1$-$3$ type hyperbolic tetrahedron and replacing the affine gluing homeomorphism by isometries preserving the decorations. Denote $E$ by the set of edges in $\CT$. We remark that the sides of the truncated triangles are not considered as the edges of the triangulation. Each edge in $E$ contains at least one ideal vertex, and is sometimes also called an \emph{internal edge} for the sake of clarity. The \emph{curvature} of the metric assigned to each internal edge $e$ is equal to $2\pi$ ($\pi$, resp.) minus the cone angle at $e$, if the gluing of tetrahedra around $e$ is cyclical (not cyclical, resp.). The \emph{length} of an internal edge refers to the directed distance between a pair of horospheres or between a horoshphere and a vertex of the truncated triangle. Our rigidity results state that for a fixed triangulation, the curvature determines the edge lengths and hence those decorated hyperbolic polyhedral metrics.
\begin{thm}\label{thm 1.2}
Let $(M,\CT)$ be a triangulated compact pseudo 3-manifold with a $1$-$3$ type triangulation $\CT$. Then a decorated $1$-$3$ type hyperbolic polyhedral metric on $(M,\CT)$ is determined up to isometry and change of decorations by its curvature.
\end{thm}
Hyperbolic polyhedral metrics have cone singularities. There had been many important works on rigidity of hyperbolic cone metrics on 3-manifolds, like the works of Hodgson-Kerckhoff \cite{Hodgson1998}, Weiss \cite{Weiss2013}, Mazzeo-Montcouquiol \cite{Mazzeo2011} and Fillastre-Izmestiev \cite{Fillastre2009, Fillastre2011,Izmestiev} and others. The difference between their rigidity and ours is that, we consider the case where the singularity consists of complete geodesics from cusp to cusp or geodesics orthogonal to the totally geodesic boundary with possible cone singularities. In addition, Luo-Yang \cite{Luo2018} had obtained the rigidity of hyperbolic polyhedral 3-manifolds produced by isometric gluing of decorated $0$-$4$ type (or $4$-$0$ type, resp.) hyperbolic tetrahedra previously. So Theorem \ref{thm 1.2} generalize Luo-Yang's rigidity to $1$-$3$ type triangulation. In the forthcoming paper \cite{Feng2023}, we will generalize Luo-Yang's rigidity to the most general setting.

Our second result concern \emph{the volume optimization program of angle structures} initiated by Casson and Rivin. Their program aims to solve Thurston's gluing equations using angle structures, which is a system of dihedral angles assigns to each edge in the triangulation (see \cite{Lackenby2000,Rivin1994}). For each angle structure, using Milnor's Lobachevsky function \cite{MR1277810}, one may assign a \emph{volume} which is the sum of volumes of all hyperbolic tetrahedra in the triangulation. Although the angle structure is weaker than hyperbolic structure, it was known for Casson and Rivin, the existence of an angle structure implies the existence of a hyperbolic structure. Moreover, if some positive angle structure maximize the volume, then there exists a geometric triangulation of the hyperbolic manifold realizing the angle structure. We refer to \cite{Futer2011,Lackenby2000,Purcell,Rivin1994,Rivin2003} for detailed proofs. Angle structures can be generalized to triangulated 3-manifolds with higher genus boundaries, see \cite{Luo2007,Luo2013,GJZ1,GJZ2} for example. Given a $1$-$3$ type triangulation $\CT$ on $M$, we say $\CT$ is \emph{closed} if the gluing of tetrahedra around each internal edge is cyclical, or equivalently, each face in $\CT$ must be pasted onto another face in $\CT$. Note each truncated triangle is not considered as a face of the triangulation $\CT$. An \emph{angle structure} corresponding to $\CT$ assigns to each internal edge in a tetrahedron a real number in $(0, \pi)$, which is called the dihedral angle of the edge in a tetrahedron, so that the sum of angles around each internal edge is $2\pi$, the sum of angles at three internal edges from each ideal vertex is $\pi$ and the sum of angles at three internal edges from each hyperideal vertex is less than $\pi$. If all the dihedral angles are allowed to be taken in the closed interval $[0, \pi]$, it is said to be a \emph{non-negative angle structure}. As a contrast, an angle structure is also called a positive angle structure in this paper. We have:

\begin{thm}\label{thm 1.3}
Suppose $(M,\CT)$ is a compact pseudo 3-manifold with a $1$-$3$ type closed triangulation $\CT$, which supports a positive angle structure. If $\alpha$ is a non-negative angle structure and $\alpha$ maximizes the volume in the space of all non-negative angle structures on $(M,\CT)$, then there exists an edge length assignment $l(e)$ to each internal edge $e\in E$, so that for each tetrahedron $\sigma$ in $\CT$,
\begin{enumerate}
\item if all angles in $\sigma$ are positive, then $\alpha$ are the dihedral angles of the decorated $1$-$3$ type hyperbolic tetrahedron of edge length given by $l$, and
\item if one angle of $\sigma$ in $\alpha$ is 0, then all angles of $\sigma$ in $\alpha$ are $0,0,0,0,\pi,\pi $ and the numbers assigned by $l$ to the edges of $\sigma$ are not the lengths of any decorated $1$-$3$ type tetrahedron.
\end{enumerate}
Conversely, if $l: E \to \IR$ is any edge length function so that (1) and (2) holds, then the corresponding angle $\alpha$ of $l$ defined by (1) and (2) maximizes the volume.
\end{thm}

Next we assume $M$ is a compact manifold with boundary, no component of which is a 2-sphere. Let $\CT$ be an ideal triangulation of $M$. The existence of an angle structure has strong restrictions for the topology of $M$. Actually, according to the pioneering work of Casson, Rivin and Lackenby (see Corollary 4.6 in \cite{Lackenby2000} and Theorem 2.2 in \cite{Lackenby2008} for example), the existence of a positive angle structure implies $M$ is irreducible, atoroidal and not Seifert fibred. Hence $M$ can be hyperbolized by Thurston's hyperbolization theorem. As a comparison, we first cite Casson-Rivin's volume optimization program (see Theorem 1.1, Theorem 1.2 and their proofs in \cite{Lackenby2008} and \cite{Futer2011}) as follows.
\begin{thm}[Casson]
\label{thm-casson}
Let $M$ be an orientable 3-manifold with boundary consisting of tori, and let $\CT$ be an ideal triangulation of $M$. Denote $\CA(\CT)$ by the set of angle structures on $(M, \CT)$. If $\CA(\CT)\neq \emptyset$, then $M$ admits a complete hyperbolic metric.
\end{thm}
In other words, the existence of an angle structure $\alpha\in\CA(\CT)$ implies the existence of a hyperbolic structure. However, it must be clarified that $\alpha$ is generally not directly related to this hyperbolic structure of $M$. If the hyperbolic structure makes the ideal triangulation $\CT$ geometric, so that the dihedral angles of those hyperbolic tetrahedra is exactly equal to the angle structure $\alpha$, then $\alpha$ is said to \emph{correspond} to the hyperbolic structure. Casson and Rivin independently proved
\begin{thm}[Casson, Rivin]
\label{thm-casson-rivin}
Let $M$ be an orientable 3-manifold with boundary consisting of tori, and let $\CT$ be an ideal triangulation of $M$. Then a point $\alpha\in\CA(\CT)$ corresponds to a complete hyperbolic metric on the interior of $M$ if and only if $\alpha$ is a critical point of the functional $vol:\CA(\CT)\rightarrow \IR$.
\end{thm}

Theorem \ref{thm 1.3} implies that Casson-Rivin's Theorem \ref{thm-casson-rivin} is still valid for $1$-$3$ type ideal triangulated 3-manifolds (part (1) of the following Theorem is known to Lackenby \cite{Lackenby2008}).
\begin{thm}\label{cor-main-1}
Suppose $M$ is an oriented compact 3-manifold with boundary, no component of which is a 2-sphere, and $\CT$ is a $1$-$3$ type ideal triangulation of $M$. Denote $\CA(\CT)$ by the set of angle structures on $(M, \CT)$.
\begin{enumerate}
  \item If $\CA(\CT)\neq \emptyset$, then $M-\partial_{t}$ admits a finite-volume complete hyperbolic structure with totally geodesic boundary, where $\partial_{t}$ is the toral boundary components of $M$.
  \item If $\alpha\in\CA(\CT)$, then $\alpha$ corresponds to a complete hyperbolic metric on $M-\partial_t$ if and only if $\alpha$ is a critical point of the volume functional.
\end{enumerate}
\end{thm}

0-efficient and 1-efficient triangulations of 3-manifolds were introduced and extensively studied in Jaco-Rubinstein \cite{Jaco}, Kang-Rubinstein \cite{Kang} and  Garoufalidis-Hodgson-Rubinstein-Segerman \cite{Garouf}, which have applications to finiteness theorems, knot theory, Dehn fillings, decision problems, algorithms, computational complexity and Heegaard splittings. By definition, an ideal triangulation $\CT$ (may not be $1$-$3$ type) of an orientable 3-manifold $M$ is \emph{0-efficient} if there are no embedded normal 2-spheres or one-sided projective planes. In addition, $\CT$ is \emph{1-efficient} if it is 0-efficient, the only embedded normal tori are vertex-linking and there are no embedded one-sided normal Klein bottles. $\CT$ is \emph{strongly 1-efficient} if there are no immersed normal 2-spheres, projective planes or Klein bottles and the only immersed normal tori are coverings of the vertex-linking tori. Similar to Theorem 1.5 in \cite{Garouf}, Theorem 2.5 and Theorem 2.6 in \cite{Kang}, we have
\begin{thm}\label{cor-main-2}
Suppose $M$ is an oriented compact 3-manifold with boundary, no component of which is a 2-sphere, and $\CT$ is a $1$-$3$ type ideal triangulation of $M$.
\begin{enumerate}
  \item If $\CT$ supports a non-negative angle structure and $M$ is atoroidal, then $\CT$ is $1$-efficient and hence $0$-efficient.
  \item If $\CT$ supports an angle structure, then $\mathcal T$ is strongly $1$-efficient and therefore $1$-efficient and $0$-efficient.
\end{enumerate}
\end{thm}

On the whole, we follow the spirit of Luo-Yang \cite{Luo2018} to prove Theorem \ref{thm 1.2} and Theorem \ref{thm 1.3}. The major differences come from the fact that the $1$-$3$ type hyperbolic tetrahedron and the $4$-$0$ (or $0$-$4$ type) tetrahedron have completely different geometry. Compared with the $0$-$4$ type hyperbolic tetrahedron, the $1$-$3$ type hyperbolic tetrahedron does not have good symmetry of dihedral angles, and compared with the $4$-$0$ type hyperbolic tetrahedron, it does not have mutual determinism between edge lengths and dihedral angles. Thus, we have to establish a new argument for the relationship between dihedral angles and decorated edges.

The paper is organized as follows. In Section 2, we study the geometry of the decorated $1$-$3$ type
hyperbolic tetrahedron. 
In Section 3, we discuss the rigidity of decorated $1$-$3$ type hyperbolic metric and prove Theorem~\ref{thm 1.2}.
In Section 4, we investigate the relationship between volume maximization and angle structures of $1$-$3$ type decorated hyperbolic polyhedral metrics
and prove Theorem~\ref{thm 1.3}.

~

\noindent
\textbf{Acknowledgements:}
The authors are grateful to Professor Gang Tian for constant encouragement. The first author would like to express special thanks to Professor Ying Zhang and Faze Zhang for many inspiring conversations and suggestions. The first author would also like to express special thanks to Professor Feng Luo, Tian Yang, Reifeng Qiu, who have always given the first author a lot of care, support, and encouragement. K. Feng is supported by NSFSC, no. 2023NSFSC1286.
H. Ge is supported by NSFC, no.12341102, no.12122119.

\section{$1$-$3$ type hyperbolic tetrahedra}
\subsection{Description of the $1$-$3$ type hyperbolic tetrahedra}
In 3-dimensional hyperbolic space ${\IH}^3$, a compact hyperbolic tetrahedron $\sigma$ is the convex hull of four non-collinear distinct points, and the four points in ${\IH}^3$ are the vertices of $\sigma$.
If these four points are all lie on $\partial{\IH}^3$, their convex hull $\sigma$ in ${\IH}^3$ is often called an ideal hyperbolic tetrahedron or \emph{ideal tetrahedron}, and each vertex is called an ideal vertex. Strictly speaking, $\sigma$ has no vertices since the four ideal vertices are not in $\sigma$. However, for the sake of convenience, we still call these four ideal vertices the vertices of $\sigma$.
Following \cite{Bao2002,Fujii1990,Schlenker2002}, a \emph{hyperideal tetrahedron} $\sigma$ in ${\IH}^3$ is a compact convex polyhedron that is diffeomorphic to a truncated tetrahedron in $\mathbb{E}^3$ with four hexagonal faces right-angled hyperbolic hexagons. The four hexagonal faces are called the \emph{faces}, and the four triangular faces isometric to hyperbolic triangles are called \emph{external faces}. An edge in a hyperideal tetrahedron is the intersection of two faces, and an \emph{external edge} is the intersection of a face and an external face. The dihedral angle at an edge
is the angle between the two (hexagonal) faces adjacent to it. The external faces are isometric to
hyperbolic triangles, and the dihedral angle between a (hexagonal) face and an external face is always $\pi/2$.

\begin{figure}[htbp]
	\centering
	\subfloat{\includegraphics[width=.8\columnwidth]{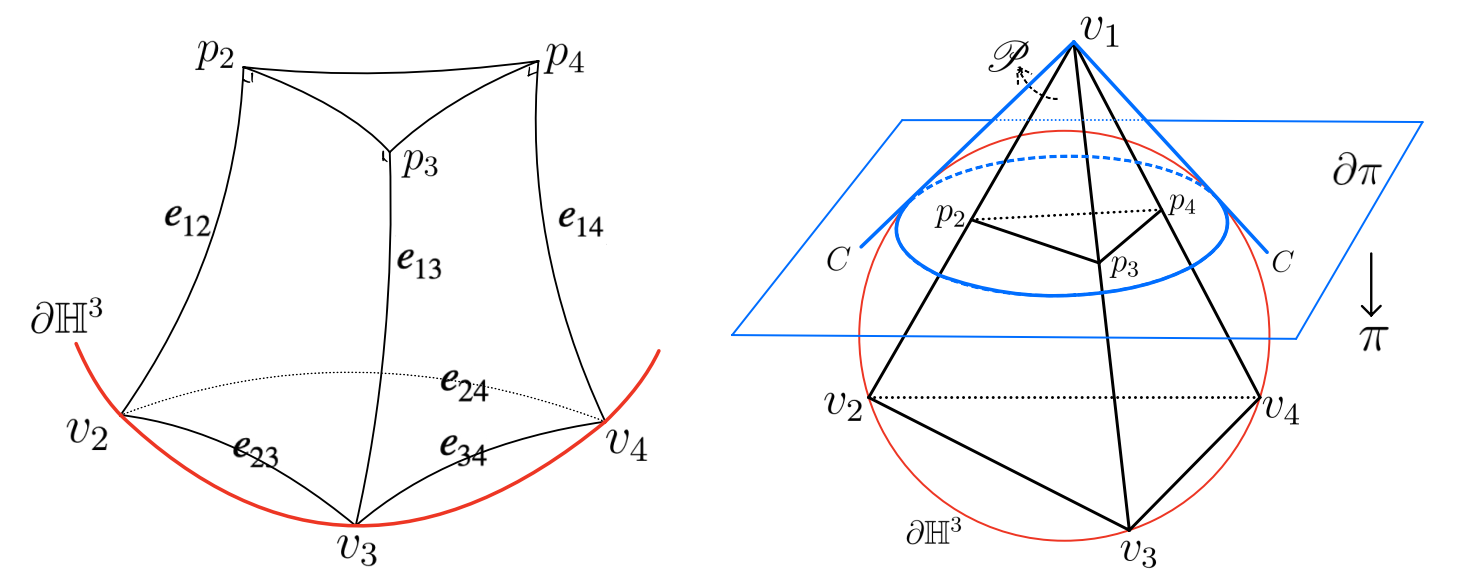}}\hspace{14pt}
	\caption{Configuration of a 1-3 type hyperbolic tetrahedron}
    \label{fig-1-3-configu}
\end{figure}

Informally, a $1$-$3$ type hyperbolic tetrahedron is a mixture of $3/4$ of an ideal tetrahedron and $1/4$ of a hyperideal tetrahedron. It  ``has" three ideal vertices and a hyperideal vertex represented by an external face derived from truncating the hyperideal vertex. However, neither these three ideal vertices, nor the hyperideal vertex, are not real entities. To be precise, we can construct a $1$-$3$ type hyperbolic tetrahedron in $\IH^3$ as follows (or see \cite{Bobenko2015}). From three vertices $p_2$, $p_3$, $p_4$ of a triangle $\triangle p_2p_3p_4$ in the hyperbolic space $\IH^3$, three rays run orthogonally to the plane of the triangle until they intersect the infinite boundary $\partial\IH^3$ in three ideal points $v_2$, $v_3$, $v_4$. See Figure \ref{fig-1-3-configu} (left). The convex hull $\sigma$ of these six points is a prism with three ideal vertices and right dihedral angles at the base triangle, it is a $1$-$3$ type hyperbolic tetrahedron. The basic triangular face is called the external face (also called a vertex triangle) of $\sigma$, and an (internal) edge in the $1$-$3$ type hyperbolic tetrahedron $\sigma$ is the geodesic ray $p_iv_i$ ($i=2, 3, 4$) or geodesic line $v_i v_j$ ($2\leq i, j\leq 4$), and an external edge (also called a vertex edge) $p_i p_j$ is the edge of the vertex triangle. By the construction, the dihedral angle between a quadrilateral face and the vertex triangle is obviously $\pi/2$.

Hyperbolic tetrahedra are best described in Klein's model for ${\IH}^3$. As described in the introduction, ${\IH}^3 \subset \IR^3\subset \IR P^3$ is the open ball representing $\IH^3$ via the Klein model. In this model, geodesics of ${\IH}^3$ correspond to the intersection of straight lines of $\mathbb{R}^3$ with ${\IH}^3$, and the totally geodesic planes in ${\IH}^3$ are the intersection of linear planes with ${\IH}^3$. Then we can obtain a $1$-$3$ type hyperbolic tetrahedron in the following way. Let $\mathscr{P}\subset \IR^3$ be a convex Euclidean tetrahedron (also considered as a projective tetrahedron) so that the vertex $v_1$ lies in $\IR^3\backslash \IH^3$ and the other vertices $v_2$, $v_3$, $v_4$ lie on $\partial\IH^3$. Suppose each edge from $v_1$ intersects $\IH^3$. Let $C$ be the cone with the apex $v_1$ tangent to $\partial\IH^3$ and $\pi$ be the half-space not containing $v_1$ so that $\partial \pi \cap \partial \IH^3= C \cap \partial \IH^3$. Then a $1$-$3$ type hyperbolic tetrahedra is given by $P := \mathscr{P}\cap \pi$, and $v_1$ is called the hyperideal vertex of $P$, and the other three vertices are called ideal vertices of $P$. Note the plane $\partial \pi$ is the dual of $v_1$ in the hyperbolic-de Sitter duality (see \cite{Schlenker2002}). Thus $P$ may be considered as derived from truncating $\mathscr{P}$ at the hyperideal vertex $v_1$ by $\partial\pi$, and the truncated face, i.e. the hyperbolic triangle $\mathscr{P}\cap \partial\pi$ is not considered as a face of $P$. $P$ has only four faces, three quadrilateral faces and one ideal triangle face. See Figure \ref{fig-1-3-configu} (right).

\begin{figure}[htbp]
	\centering
	\subfloat[]{\includegraphics[width=.40\columnwidth]{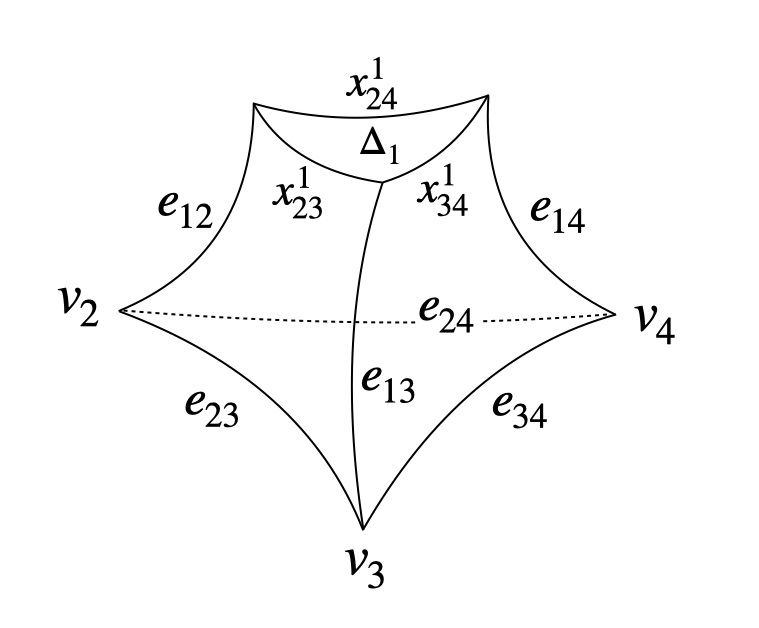}}\hspace{5pt}
	\subfloat[]{\includegraphics[width=.40\columnwidth]{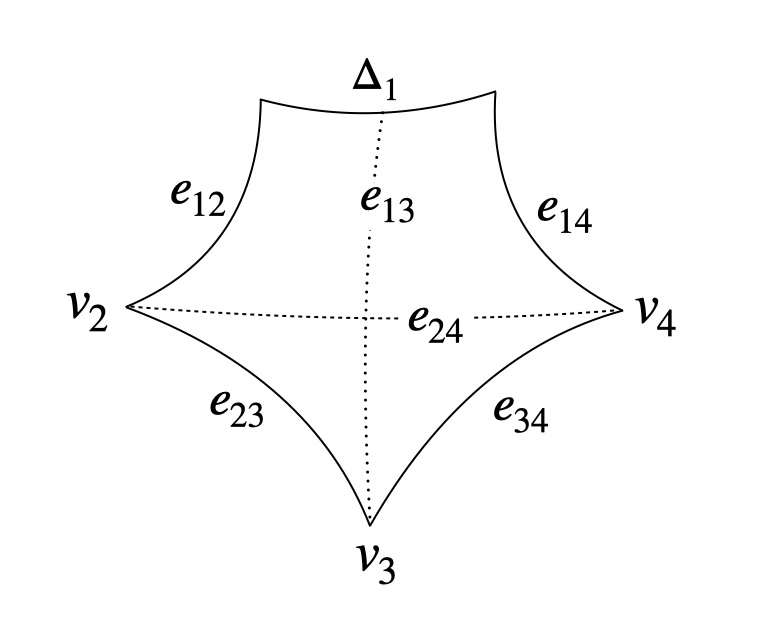}}
	\caption{left (right): a $1$-$3$ type (flat $1$-$3$ type) hyperbolic tetrahedron}
    \label{tetra-flat-tera}
\end{figure}

Now, we introduce some notations in a $1$-$3$ type hyperbolic tetrahedra. For example, in Figure \ref{tetra-flat-tera} (A), $\sigma$ is a $1$-$3$ type hyperbolic tetrahedron in ${\IH}^3$, without loss of generality, we can always assume that the vertex $v_1$ is the hyperideal vertex, and $v_2$, $v_3$, $v_4$ are ideal vertices. We denote the vertex triangle by $\triangle_1$. For $\{i,j,k\} = \{2,3,4\}$, denote $e_{1i}$, $e_{1j}$, $e_{1k}$ by the edge joining $\triangle_1$ to $v_i$, $v_j$, $v_k$, denote $e_{jk}$ by the edge joining $v_j$ and $v_k$, and denote $Q_{1jk}$ by the quadrilateral face adjacent to $e_{1j}$, $e_{1k}$, $e_{jk}$. The vertex edge $\triangle_1 \cap Q_{1jk}$ is denoted
by $x^1_{jk}$. For $p, q\in\{1, 2,3,4\}$, $p\neq q$, the dihedral angle at the internal edge $e_{pq} $ is denoted by $\alpha_{pq}$. As a convention, we always assume $\alpha_{pq}= \alpha_{qp}$.

To prove the main result, we must face the phenomenon of degradation. At this point, the vertex triangle $\triangle_1=\triangle p_2p_3p_4$ degenerates to a geodesic segment, and we get a \emph{flat} $1$-$3$ type hyperbolic tetrahedron, see Figure \ref{tetra-flat-tera} (B) for instance, where the vertex $p_3$ falls between the geodesic segment $p_2p_4$.

\subsection{Geometry in the $1$-$3$ type hyperbolic tetrahedra}
Obviously, the hyperbolic length of each internal edge is infinity. In order to describe the geometry in a hyperbolic tetrahedron with ideal vertices, i.e. the relationship between hyperbolic edge lengthes and dihedral angles, infinite edge length is inconvenient. A very nice approach is to adopt Penner's strategy \cite{Penner1987}, that is attaching a horosphere (called a decoration) at each ideal vertex. By definition, a decorated $1$-$3$ type hyperbolic tetrahedron is a pair of $(\sigma, \{H_2, H_3, H_4\})$, where $\sigma$ is a $1$-$3$ type hyperbolic tetrahedron and $H_i$ is a horosphere centered at the ideal vertex $v_i$ for $i\in\{2,3,4\}$.
We call $\{H_2, H_3, H_4\}$ the decorations of $\sigma$. For a pair of ideal vertices $v_i$ and $v_j$ with $i, j\in\{2,3,4\}$, $i\neq j$, the signed edge length $l_{ij}$ of the decorated edge $e_{ij}$ is defined as follows. The absolute value $\vert l_{ij} \vert $ is the distance between $H_i\cap e_{ij}$ and $H_j\cap e_{ij}$, so that $l_{ij}>0$ if $H_i$ and $H_j$ are disjoint and $l_{ij}\leq 0$ if $H_i\cap H_j\neq\emptyset$. For each ideal vertex $i$ with $i\in\{2,3,4\}$, the sighed edge length $l_{1i}$ is the hyperbolic length of the edges $e_{1i}$ between $p_i$ and $H_i$. It is a positive number if and only if $p_i$ falls outside the horosphere $H_i$. 

\begin{figure}[htbp]
	\centering
	\subfloat[$l_{ij}>0$]{\includegraphics[width=.32\columnwidth]{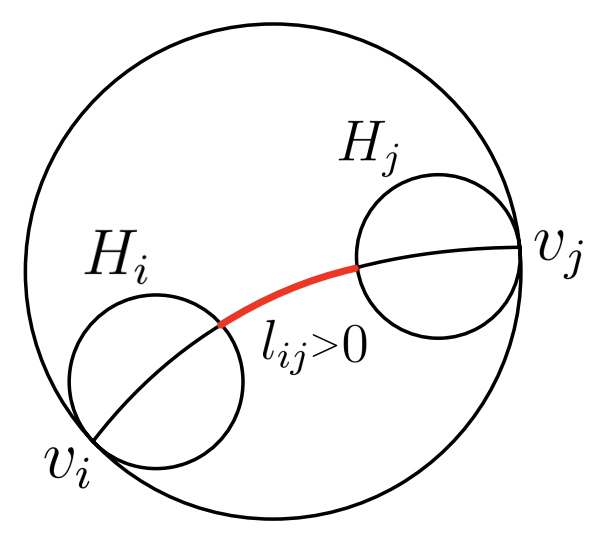}}\hspace{14pt}
	\subfloat[$l_{ij}<0$]{\includegraphics[width=.32\columnwidth]{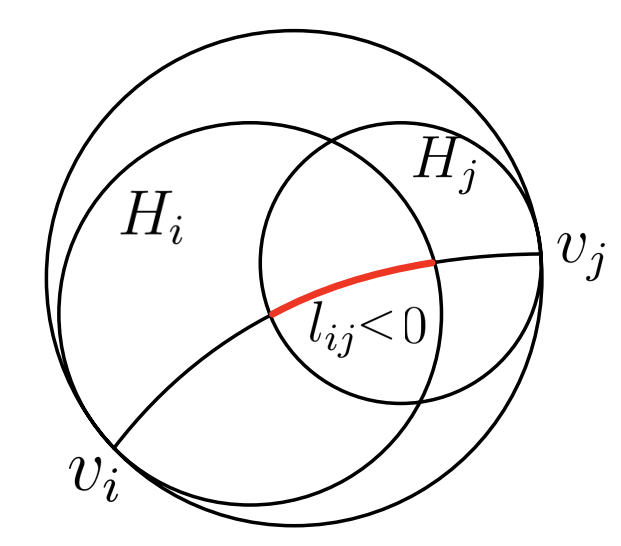}}
	\caption{decorated edge lengths in Poincar\'{e}'s disk model }
\end{figure}


A $1$-$3$ type hyperbolic tetrahedron is uniquely determined by its six dihedral angles, and a decorated $1$-$3$ type hyperbolic tetrahedron is uniquely determined by its six edge lengthes. In fact, for the moduli space of the dihedral angles and the moduli space of the decorated edge lengths of a decorated $1$-$3$ type hyperbolic tetrahedron, we have:
\begin{prop}[\cite{Bao2002,Fujii1990,Schlenker2002}]\label{prop 2.1}
Let $\sigma$ be a $1$-$3$ type hyperbolic tetrahedron in $\IH^3$, and $v_1 \in \sigma$ is the hyperideal vertex, then the isometry class of $\sigma$ is determined by its six dihedral angles $(\alpha_{12},\ldots, \alpha_{34})\in{\IR}^6_{>0}$ which satisfies the condition that (recall $\alpha_{ij}=\alpha_{ji}$)
\begin{equation}
\label{condition-angle}
  \left\{
    \begin{aligned}
      \alpha_{12} + \alpha_{13}+ \alpha_{14} < \pi \\
      \alpha_{21}+ \alpha_{23}+ \alpha_{24} =\pi \\
      \alpha_{31}+ \alpha_{32}+ \alpha_{34} =\pi \\
      \alpha_{41}+ \alpha_{42}+ \alpha_{43} =\pi
    \end{aligned}
  \right.
\end{equation}
Conversely, given $(\alpha_{12},\ldots,\alpha_{34})\in{\IR}^6_{>0}$ which satisfies the above condition (\ref{condition-angle}), then there exists a $1$-$3$ type hyperbolic tetrahedron having $\alpha_{ij}$ as its dihedral angles.

If further assume $\sigma$ is decorated, then the isometry class of $\sigma$ is determined by its edge length vector $l=(l_{12},\ldots,l_{34})$ up to change of decorations.
\end{prop}

Due to Proposition \ref{prop 2.1}, the space of isometry class of $1$-$3$ type hyperbolic tetrahedron parametrized by dihedral angles is an open convex polytope in ${\IR}^6_{>0}$:
\begin{equation}
\CB=\Bigl\{(\alpha_{12},\ldots, \alpha_{34})\in{\IR}^6_{>0}:\sum_{j\neq 1}\alpha_{1j}<\pi,\;\sum_{j\neq i} \alpha_{ij}=\pi,\;\forall\;i \in \{2,3,4\}\Bigr\}.
\end{equation}

Using the condition (\ref{condition-angle}), we obtain the following equalities of dihedral angles:
\begin{equation} \label{angle}
  \left\{
    \begin{aligned}
      \alpha_{23}  =  -\frac{1}{2}( \alpha_{12}+ \alpha_{13}- \alpha_{14} ) + \frac{\pi}{2} \\
      \alpha_{24}  =  -\frac{1}{2}( \alpha_{12}+ \alpha_{14}- \alpha_{13} ) + \frac{\pi}{2} \\
      \alpha_{34}  =  -\frac{1}{2}( \alpha_{13}+ \alpha_{14}- \alpha_{12} ) + \frac{\pi}{2}
    \end{aligned}
  \right.
\end{equation}
that is, for $\{i,j,k\} = \{2,3,4\}$,
\begin{equation}
  \alpha_{ij}  =  -\frac{1}{2}( \alpha_{1i}+ \alpha_{1j}- \alpha_{1k} ) + \frac{\pi}{2} .
\end{equation}
As a consequence,
\[ \alpha_{23} + \alpha_{24} + \alpha_{34} = \frac{1}{2} (3\pi - ( \alpha_{12}+ \alpha_{13}+ \alpha_{14})) , \]
hence $\pi<\alpha_{23}+\alpha_{24}+\alpha_{34}<\dfrac{3\pi}{2}$.\\



Next, we concern the cosine and sine laws in a decorated $1$-$3$ type hyperbolic tetrahedron. Let $(l_{12},\ldots,l_{34} )\in \IR^6 $ be the decorated edge length vector of a decorated $1$-$3$ type hyperbolic tetrahedron $\sigma$. For $\{i, j, k\}=\{2,3,4\}$, denote $\theta_{ij}^1$ by the hyperbolic length of the vertex edge $x^1_{ij}$. By the cosine law in the hyperbolic triangle $\triangle_1=\triangle p_2p_3p_4$, we have
\begin{equation}
\label{solv-triangle-P}
  \cos \alpha_{1i} = \frac{\cosh \theta_{ij}^1\; \cosh \theta_{ik}^1 - \cosh \theta_{jk}^1}{\sinh \theta_{ij}^1 \; \sinh\theta_{ik}^1}.
\end{equation}
Further by the cosine and sine formulas in the quadrilateral $Q_{1ij}$ (for instance, see \cite{Bobenko2015,Leibon2002}, or the Appendix A in Guo-Luo \cite{Guo2009}), we have:
\begin{equation}
\label{Q1ij}
\sinh^2 \, \frac{\theta_{ij}^1}{2} = e^{l_{ij}- l_{1i} - l_{1j}} .
\end{equation}
Thus using the above (\ref{Q1ij}), we obtain
\begin{equation}
\label{cosh-xita}
  \left\{
    \begin{array}{ll}
      \cosh  \theta_{ij}^1= 2e^{l_{ij}-l_{1i}-l_{1j}}+1\\[10pt]
      \sinh  \theta_{ij}^1= 2\sqrt{e^{l_{ij}-l_{1i}-l_{1j}}}\cdot \sqrt{1+e^{l_{ij}-l_{1i}-l_{1j}}} .
    \end{array}
  \right.
\end{equation}

\begin{figure}[htbp]
	\centering
	\subfloat{\includegraphics[width=.98\columnwidth]{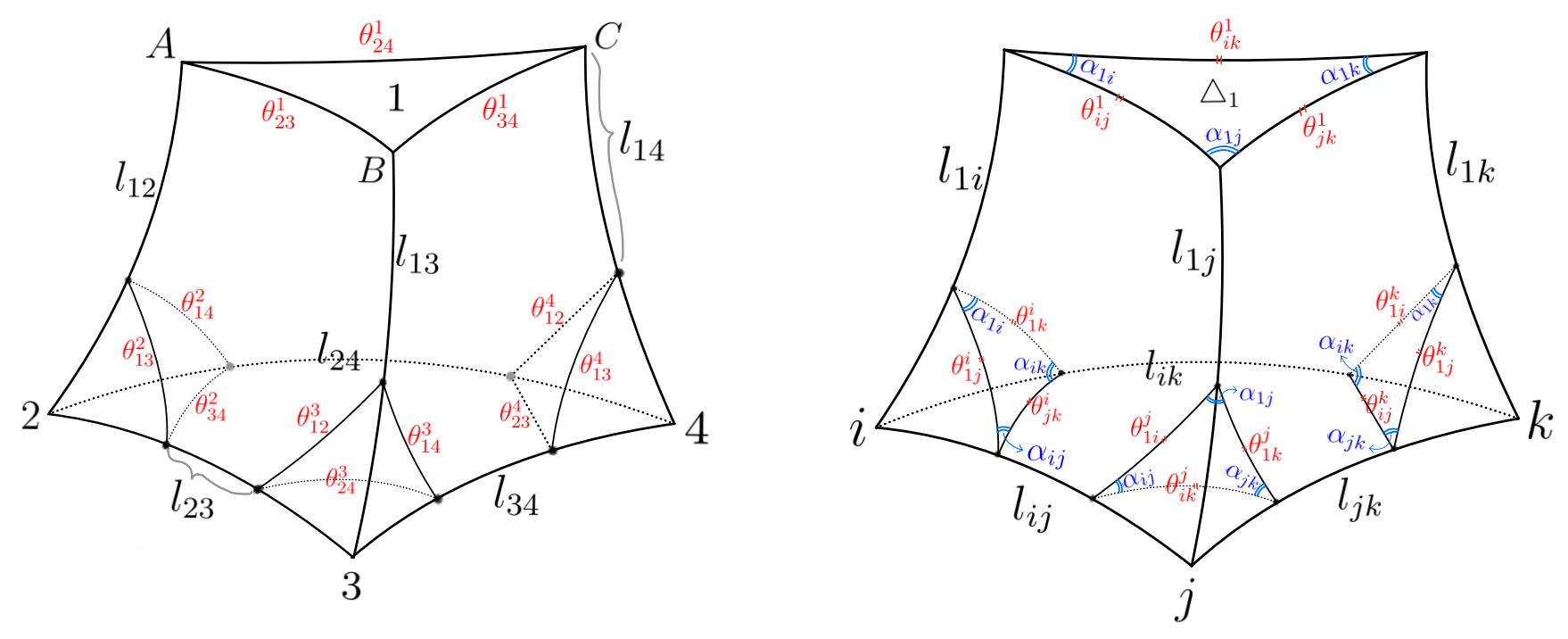}}\hspace{14pt}
	\caption{Notations in a $M_{1|3}$ tetrahedron}
    \label{fig-1-3-label}
\end{figure}
As above, for $\{i, j, k\}=\{2,3,4\}$, recall $H_i$ is a horoshpere centered at an ideal vertex $v_i$. Consider the intersection $H_i\cap\sigma$, it is an Euclidean triangle, with three edges denoted by $x_{jk}^i$, $x_{1j}^i$ and $x_{1k}^i$. Here $x_{jk}^i$ is the Euclidean arc in $H_i\cap\sigma$ connecting the edges $e_{ij}$ and $e_{ik}$. Similarly, $x_{1j}^i$ ($x_{1k}^i$ resp.) is the Euclidean arc connecting $e_{1i}$ and $e_{ij}$ ($e_{ik}$ resp.). See Figure \ref{fig-1-3-label} (right) for illustration. Denote $\theta_{jk}^i$ by the Euclidean length of $x_{jk}^i$ in $H_i\cap\sigma$, and denote $\theta_{1j}^i$ and $\theta_{1k}^i$ in the same way.
Using the Appendix A in Guo-Luo \cite{Guo2009}, we have
\begin{equation}
\label{triangle-XYZ}
  \left\{
    \begin{array}{ll}
      \theta_{1j}^i=2\sqrt{(e^{l_{ij}-l_{1i}}+e^{l_{1j}})/e^{l_{ij}+l_{1i}}}\\[10pt]
      \theta_{jk}^i=2\sqrt{e^{l_{jk}-l_{ij}-l_{ik}}}
    \end{array}
  \right.
\end{equation}

\begin{lemma}\label{angleedge}
Let the vector $l = (l_{12},l_{13},l_{14},l_{23},l_{24},l_{34}) \in \IR^6$ be a independent variable, define the functions $\phi_{pq}(l)$ for $p, q\in\{1, 2,3,4\}$, $p\neq q$ as follows: set $\phi_{1i}(l)= \phi_{i1}(l)$ and
\begin{equation}
\label{phi-1i}
  \phi_{1i}(l)=
  \frac{2 e^{l_{ij}-l_{1i}-l_{1j}} e^{l_{ik}-l_{1i}-l_{1k}} + e^{l_{ij}-l_{1i}-l_{1j}} + e^{l_{ik}-l_{1i}-l_{1k}} - e^{l_{jk}-l_{1j}-l_{1k}}}{2 \sqrt{e^{l_{ij}-l_{1i}-l_{1j}}}\cdot \sqrt{1+e^{l_{ij}-l_{1i}-l_{1j}}}\cdot \sqrt{e^{l_{ik}-l_{1i}-l_{1k}}}\cdot \sqrt{1+e^{l_{ik}-l_{1i}-l_{1k}}}}
\end{equation}
for $i \in \{2,3,4\}$, and further define
\begin{equation}
\label{phi-jk}
\phi_{jk}(l)=\frac{e^{l_{ij}+l_{1k}}+e^{l_{ik}+l_{1j}}-e^{l_{jk}+l_{1i}}}{2\sqrt{e^{l_{ij}+l_{ik}+l_{1j}+l_{1k}}+e^{l_{jk}+l_{ij}+l_{ik}}}}
\end{equation}
for $\{i,j,k\}=\{2,3,4\}$. Then $\phi_{pq}(l)=\phi_{qp}(l)$ for each pair of $p, q\in\{1, 2,3,4\}$, $p\neq q$. Moreover, if the vector $l$ is the decorated edge length vector of some decorated $1$-$3$ type hyperbolic tetrahedron $(\sigma, \{H_2, H_3, H_4\})$, then $\alpha_{pq}=\alpha_{pq}(l)$ is a function of $l$ and
\begin{equation}
\phi_{pq}(l)=\cos({\alpha_{pq}}(l)).
\end{equation}
\end{lemma}
\begin{proof}
For $i \in \{2,3,4\}$, $\phi_{1i}(l)=\phi_{i1}(l)$ by the settings. By (\ref{phi-jk}), $\phi_{jk}(l)$ is symmetric with respect to the lower indicators $j$ and $k$ for $\{i,j,k\}=\{2,3,4\}$. Overall, $\phi_{pq}(l)=\phi_{qp}(l)$ for each pair of $p, q\in\{1, 2,3,4\}$, $p\neq q$.

In the Euclidean triangle $H_i\cap\sigma$, using the cosine law we get (see Figure \ref{fig-1-3-label} (right))
\begin{equation}
\label{alpha-ij}
\cos{\alpha_{ij}}=\frac{{(\theta_{1j}^{i})}^{2}+{(\theta_{jk}^i)}^{2}-{(\theta_{1k}^i)}^{2}}{2\theta_{1j}^i\cdot\theta_{jk}^i}.
\end{equation}

Using (\ref{solv-triangle-P})-(\ref{triangle-XYZ}) and (\ref{alpha-ij}), $\alpha_{pq}(l)$ is determined by $l$ and hence a function of $l$. Moreover, by direct calculations, the equalities $\phi_{pq}(l)=\cos({\alpha_{pq}}(l))$ could be checked by substituting (\ref{cosh-xita}) into (\ref{solv-triangle-P}) and substituting (\ref{triangle-XYZ}) into (\ref{alpha-ij}).
\end{proof}

\begin{rmk}
\label{domain-of-phi}
By (\ref{phi-1i})-(\ref{phi-jk}), the domain of $\phi_{pq}(l)$ is the whole space $\IR^6$. However, the domain of $\alpha_{pq}(l)$ is only a proper subset of $\IR^6$ (see Proposition \ref{prop 4.4}), hence $\phi_{pq}$ may be considered as an extension of the function $\cos({\alpha_{pq}}(l))$, see Section \ref{extend-angle}. This fact will play a central role in our paper.
\end{rmk}

\begin{rmk}
\label{unchange}
Given a decorated $1$-$3$ type hyperbolic tetrahedron $(\sigma, \{H_2, H_3, H_4\})$, if we change the sizes of these $H_i$, the decorated edge length will also change accordingly. However, by the formula (\ref{phi-1i})-(\ref{phi-jk}), the value of the functions $\phi_{pq}$ will not change, because the differences $\Delta l_{ij}$ of $l_{ij}$ always satisfy $\Delta l_{ij}=\Delta l_{1i}+\Delta l_{1j}$ for $ \{i,j\} \subseteq \{2, 3, 4\}$.
\end{rmk}

\subsection{The volume and its dual: co-volume}
For an ideal tetrahedron, its opposite sides have the same dihedral angle. Assuming that the dihedral angles on its three sets of opposite sides are $\theta_1$, $\theta_2$ and $\theta_3$, respectively, then its volume is $\Lambda(\theta_1)+\Lambda(\theta_2)+\Lambda(\theta_3)$, where $\Lambda(\theta)$ is Milnor's Lobachevsky function defined by
\[\Lambda(\theta) = -\int_{0}^{\theta} \ln \vert 2\sin t \vert \di t.\]
which is $\pi$-periodic, odd and smooth except at $\theta\in \pi\mathbb{Z}$, with tangents vertical (\cite{MR1277810}). Note $\theta_1+\theta_2+\theta_3=\pi$ in an ideal tetrahedron, and Rivin \cite{Rivin1994} proved that the volume is strictly concave down on $\{(\theta_1, \theta_2, \theta_3)\in(0, \pi)^3:\theta_1+\theta_2+\theta_3=\pi\}$.

As above, $(\sigma, \{H_2, H_3, H_4\})$ is a decorated $1$-$3$ type tetrahedron, and the dihedral angle at the edge $e_{pq}$ is $\alpha_{pq}$ for $\{p,q\}\subseteq \{1,2,3,4\}$. Obviously, these decorations do not affect the volume of the tetrahedron $\sigma$, which is a finite and positive. By hyperbolic doubling $\sigma$ along the vertex triangle $\triangle_1$, one get an ideal triangular prism which can be decomposed into three ideal tetrahedra. The dihedral angles in these three ideal tetrahedra can be derived from the dihedral angles of $\sigma$. By this way, one can express
the volume of $\sigma$ as a function of its six dihedral angles (see \cite{Leibon2002} or \cite{GY} for instance):
\begin{equation}
\label{volume-alpha}
2 vol(\sigma(\alpha_{12},\ldots, \alpha_{34}))=\Lambda\Bigl(\frac{\pi-\alpha_{12}-\alpha_{13}-\alpha_{14}}{2}\Bigr)+\sum_{1\leq i<j\leq 4}\Lambda(\alpha_{ij}).
\end{equation}
Consider the volume function $vol: \CB\rightarrow \mathbb{R}$, there are only three free variables by (\ref{angle}), for example $\alpha_{12}$, $\alpha_{13}$ and $\alpha_{14}$ in the six dihedral angles of $\sigma$. Similar with Rivin's observation, the volume $vol(\alpha_{12},\ldots, \alpha_{34})$ of a $1$-$3$ type hyperbolic tetrahedron is also concave down in $\CB$. If we project $\CB$ to a convex polytope $\CB_{proj}\subset\mathbb{R}^3$ along any three free variables, then the volume function $vol(\cdot)$ is strictly concave down on $\CB_{proj}$ (we refer to Section 3 of Schlenker \cite{Schlenker2002} for a more precise and comprehensive formulation). This is because:
\begin{lemma}[\cite{Leibon2002,Schlenker2002,Rivin1994}]
\label{hessian-neg-definite}
The Hessian matrix of $vol(\cdot)$ with respect to the six dihedral angles $\{\alpha_{pq}\}$ is semi-negative definite, and is negative definite with respect to any three free variables in $\{\alpha_{pq}\}$.
\end{lemma}

Let $\CL\subset\IR^6$ be the set of vectors $(l_{12}, \ldots, l_{34})$ such that there exists a decorated
$1$-$3$ type hyperbolic tetrahedron $(\sigma, \{H_2, H_3, H_4\})$ with six decorated edge lengthes $\{l_{pq}\}_{1\leq p,q\leq 4}$.
We will see that $\CL\subset \IR^6$ is a simply connected open set in the future. To see $\CL$ is not the whole space $\IR^6$, it is easy to find out that $(1,1,1,2,2,2)\in\CL$ but $(1,1,1,6,4,2+\ln c)\notin\CL$, where $c=2e^6+e^2+e^4-2e^3\sqrt{(1+e^2)(1+e^4)}$. Actually, we have:

\begin{prop}\label{prop 4.4}
The space of decorated $1$-$3$ type hyperbolic tetrahedra $(\sigma, \{H_2, H_3, H_4\})$ parameterized by the decorated edge lengths is
\begin{equation}
  \CL = \big\{(l_{12},\ldots,l_{34}) \in {\IR}^6 \mid \phi(l_{ij})\in(-1,1) \; \text{for all } \; \{i,j\} \subseteq \{1,2,3,4\}\big\} .
\end{equation}
\end{prop}

\begin{proof}
If $l$ is the edge length vector of a decorated $1$-$3$ type hyperbolic tetrahedron, then by the definition, each of the dihedral angles $\alpha_{ij}(l) \in (0, \pi)$. We apply the cosine law to the vertex triangle $\triangle_1 $ and polygons $H_{123}, \, H_{124}, \, H_{134}$, $H_{234}$, then we have $\phi(l_{ij}) = \cos(\alpha_{ij}) \in (-1,1)$ for all
$\{i,j\} \subseteq \{1,2,3,4\}$.

Conversely, for each $l \in {\IR}^6$, if we require that $\phi(l_{ij}) \in (-1,1) $, then by the definition of $\phi(l_{1i})$, ${\theta}_{ij}^1,\; {\theta}_{ik}^1,\;{\theta}_{jk}^1$ satisfy the triangular inequality. Then there exists a unique vertex hyperbolic triangle $\triangle_1 $ having them as edge lengths.
Taking ${\alpha}_{1i} = {\cos}^{-1} \phi(l_{1i}) \in (0,\pi)$, we see that $\alpha_{1i},\, \alpha_{1j},\, \alpha_{1k}$ are the inner angles of  triangle $\triangle_1 $, hence they satisfy $ \alpha_{1i}+ \alpha_{1j}+ \alpha_{1k} < \pi $. Since we have the identity that  $\alpha_{23},\, \alpha_{23},\, \alpha_{34}$ can be express linearly by
$\alpha_{12},\, \alpha_{13},\, \alpha_{14}$, then all the six dihedral angles $ ( \alpha_{12}, \ldots, \alpha_{34}) $ are determined, and by directly calculating, their cosine function value also equal with the values given by Proposition~\ref{prop 2.1}, we see that $l$ is the edge length of a $1$-$3$ type tetrahedron.
\end{proof}
Consider the volume function $vol: \CB\rightarrow \mathbb{R}$, the Schl\"{a}fli formula \cite{Bonahon1998,Schlenker2002} says
\begin{equation}
\di vol=-\frac{1}{2}\sum_{1\leq i<j\leq 4}l_{ij}d\alpha_{ij}.
\end{equation}
\begin{lemma}[co-volume]
\label{p-cov}
The co-vlume function $cov:\CL \rightarrow \IR$ defined by
\begin{equation}
cov(l_{12},\ldots,l_{34})=2vol(\alpha_{12},\ldots, \alpha_{34})+\sum_{1\leq i<j\leq 4} \alpha_{ij}\cdot l_{ij}
\end{equation}
is $C^{\infty}$-smooth and satisfies $\nabla_l cov=\alpha$, or equivalently, for each $\{i,j\}\subseteq \{1,2,3,4\}$,
\begin{equation}
\label{partial-cov}
\frac{\partial cov}{\partial l_{ij}}=\alpha_{ij}.
\end{equation}
\end{lemma}

\begin{proof}
The formula (\ref{partial-cov}) can be proved by direct computations, see the Appendix. Here we prove it by Schl\"{a}fli's formula.
Since $\di vol = -\frac{1}{2}\sum_{1 \le i < j \le 4} l_{ij}\di {\alpha}_{ij}$,
then
\[\di cov = 2 \di vol + \sum_{1 \le i < j \le 4}\big(\alpha_{ij} \di l_{ij} + l_{ij}\di \alpha_{ij}\big)= \sum_{1 \le i < j \le 4} \alpha_{ij} \di l_{ij},\]
which implies $\partial cov/\partial l_{ij}=\alpha_{ij}$ for each $\{i,j\}\subseteq \{1,2,3,4\}$.
\end{proof}

Denote $ v_1=(1,0,0,1,1,0)$, $v_2=(0,1,0,1,0,1)$ and $v_3=(0,0,1,0,1,1)$, then $\CL$, $\alpha$, $vol$ and $cov$ are all invariant along the direction $v_1$, $v_2$ and $v_3$ by Lemma~\ref{angleedge} and Remark \ref{unchange}. To be precise, denote $L(v_1,v_2,v_3)=\text{span}\{v_1, v_2, v_3 \}$, and let $l\in\CL$ be the decorated edge length vector of a decorated $1$-$3$ type hyperbolic tetrahedron $(\sigma, \{H_2, H_3, H_4\})$, then for any $a$, $b$, $c\in \IR$, the decorated edge length vector $l'=l+av_1+bv_2+cv_3$ is still in $\CL$. It represents the same $1$-$3$ type hyperbolic tetrahedron $\sigma$ but with different decorations $\{H'_2, H'_3, H'_4\}$. Except for different decorations, the geometry of $l$ and $l'$ are completely identical, and therefore $\alpha(l)=\alpha(l')$ and $vol(l)=vol(l')$. We further remark that both $\CB$ and $\CL / L(v_1,v_2,v_3)$ can be considered as simply connected open sets in $\mathbb{R}^3$, and they are diffeomorphic by Proposition \ref{prop 2.1}. Consequently, $\CL\subset \IR^6$ is simply connected.

\begin{thm}
The co-volume function $cov:\CL \rightarrow \IR$ is locally convex. To be precise, it has a semi-positive definite Hessian matrix with rank 3 at each $l\in \CL$. Moreover, $cov$ is locally strictly convex on the open set $\CL / L(v_1,v_2,v_3)\subset \mathbb{R}^3$.
\end{thm}
\begin{proof}
Consider the co-volume function $cov(l)$, by Lemma \ref{p-cov}, $\partial cov/\partial l_{ij}=\alpha_{ij}$ for each $\{i,j\}\subseteq \{1,2,3,4\}$, hence the Hessian of $cov$ with respect to $l$ is exactly ${[\partial  \alpha_{ij} / \partial l_{pq}]}_{6\times 6} $. Consider $\CB$ as an open set in $\mathbb{R}^3$ with three coordinates $\alpha_{12}$, $\alpha_{13}$ and $\alpha_{14}$, then the smooth function $vol(\alpha_{12},\alpha_{13},\alpha_{14}):\CB\rightarrow\mathbb{R}$ has a negative definite Hessian matrix at each point. We could well assume that $l_{12}$, $l_{13}$ and $l_{14}$ are the three coordinates of $\CL/L(v_1,v_2,v_3)$, and consider the quotient co-volume function $cov:\CL/L(v_1,v_2,v_3)\rightarrow\mathbb{R}$. It is smooth, and its gradient
$(\partial cov/\partial l_{12}$, $\partial cov/\partial l_{13}$, $\partial cov/\partial l_{14}):\CL/L(v_1,v_2,v_3)\rightarrow\mathbb{R}^3$
is injective by Proposition \ref{prop 2.1}. The Legendre transform (differs by a coefficient of 2) is exactly the above $vol(\alpha_{12},\alpha_{13},\alpha_{14})$. Therefor, the Hessian matrix of the quotient $cov$ (with respect to $l_{12}$, $l_{13}$ and $l_{14}$) is in fact the inverse of the Hessian matrix of $vol$ (with respect to $\alpha_{12}$, $\alpha_{13}$ and $\alpha_{14}$), but differs by a coefficient of $-2$. According to Lemma \ref{hessian-neg-definite}, the latter is negatively definite, hence the Hessian matrix of the quotient $cov$ is positive definite. As a consequence, $cov$ is locally strictly convex on the open set $\CL / L(v_1,v_2,v_3)\subset \mathbb{R}^3$. Further note $cov:\CL\rightarrow\mathbb{R}$ is linear along the space $L(v_1,v_2,v_3)$, hence $cov:\CL\rightarrow\mathbb{R}$ is locally convex in $\CL$. Since ${[\partial\alpha_{ij}/\partial l_{pq}]}_{6\times6}$ can be transformed through a congruence in to a matrix with positive definite $3\times3$ principal minor (which is the Hessian matrix of the quotient $cov$) and other entries zero, ${[\partial\alpha_{ij}/\partial l_{pq}]}_{6\times6}$ is semi-positive definite with rank 3.
\end{proof}

\begin{cor}
\label{cor-cov-a*l}
For an given $\theta\in\CB$, the following locally convex function
\begin{equation}
\psi_{\theta} (l)=cov(l)-\theta\cdot l
\end{equation}
satisfies $\psi_{\theta} (l+rv_1+sv_2+tv_3)=\psi_{\theta}(l)$
for all $r, s, t \in \IR$ and $l\in\CL$. Moreover, an decorated edge length vector $l\in\CL$ with $\alpha(l)=\theta$ is the unique (up to decorations) critical point of $\psi_{\theta}$.
\end{cor}

\subsection{Degeneration of the $1$-$3$ tetrahedra}
If $l\in\CL$ is the decorated edge length vector of a decorated $1$-$3$ type hyperbolic tetrahedron $(\sigma, \{H_2, H_3, H_4\})$, the geometric meaning of those $l_{pq}$, $\theta^r_{pq}$ and $\alpha_{pq}$ ($\{p,q,r\}\subset\{1,2,3,4\}$) are clear. How about the general $l\in\mathbb{R}^6$?

\begin{figure}[htbp]
	\centering
	\subfloat{\includegraphics[width=.75\columnwidth]{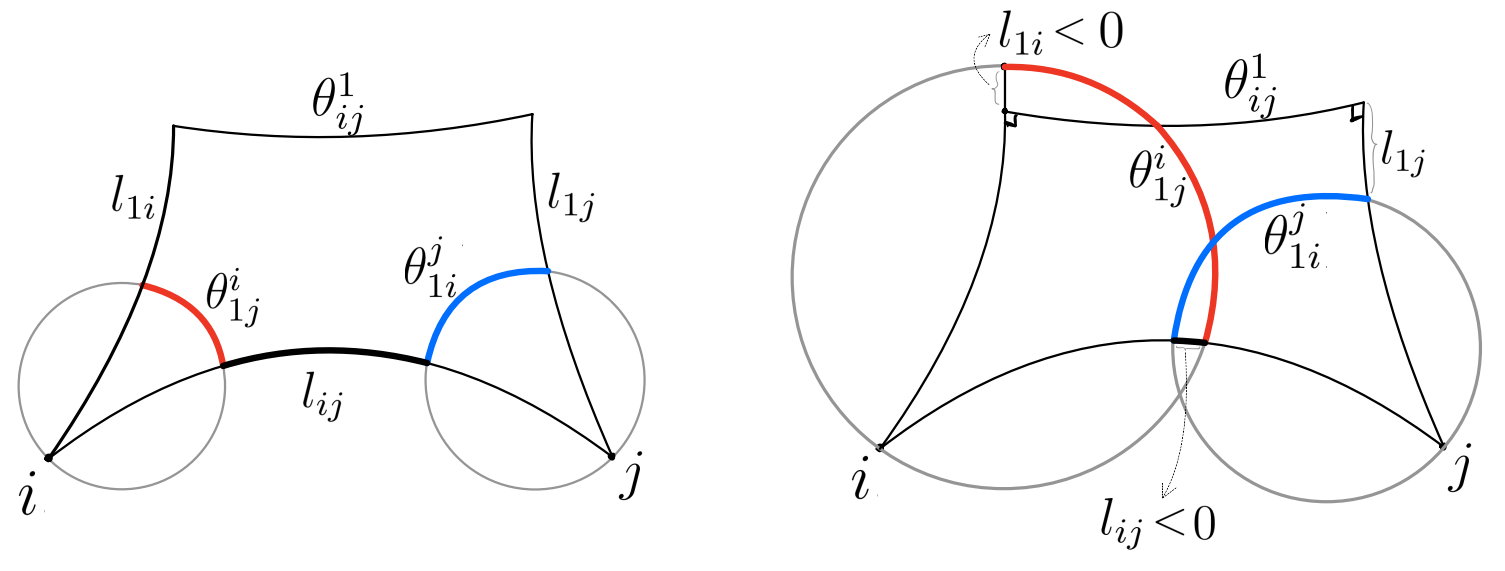}}\hspace{14pt}
	\caption{Configuration of a quadrilateral}
    \label{fig-quad-config}
\end{figure}

If $l\notin\CL$, then it represents a degenerate geometry in $\sigma$. Let's elaborate in detail. Firstly, $l$ does not have one-dimensional degeneracy, i.e. any given number $l_{pq}\in\mathbb{R}$ can always be implemented as the signed edge length of some truncated internal edge $e_{pq}$. Secondly, $l$ does not have two-dimensional degeneracy. This is because, on the one hand, for a quadrilateral $Q_{1ij}$ with $\{i,j\}\subset\{2,3,4\}$, any three real numbers $l_{ij}$, $l_{1i}$ and $l_{1j}$ can be appropriately adjusted by modifying the sizes of horospheres $H_i$ and $H_j$, so that the three numbers become positive. For such three positive real numbers, there always exists a hyperbolic quadrilateral, as shown in Figure \ref{fig-quad-config} (left), which has two right angles and two ideal vertices, each with a horosphere, so that the truncated edge lengths of the three internal edges are those three positive real numbers $l_{ij}$, $l_{1i}$ and $l_{1j}$, respectively. By adjusting the size of these two horospheres in reverse, we obtain a decorated hyperbolic quadrilateral, with the original three real numbers $l_{ij}$, $l_{1i}$ and $l_{1j}$ as signed edge lengthes. See Figure \ref{fig-quad-config} (right) for instance, where $l_{1j}>0$, $l_{1i}$ and $l_{ij}<0$.
On the other hand, for the bottom triangle $\triangle 234$, given three real numbers $l_{23}$, $l_{24}$ and $l_{34}$, there is always a unique decorated ideal triangle with three signed edge lengths of $l_{23}$, $l_{24}$ and $l_{34}$. These two aspects indicate that all four two-dimensional faces of $\sigma$ can be geometrically realized, and $l$ does not undergo two-dimensional degeneration. As a consequence, the formulas (\ref{Q1ij}), (\ref{cosh-xita}) and (\ref{triangle-XYZ}) are always meaningful for general $l$. Finally, $l$ indeed undergoes degeneration in three dimension. In fact, let's consider the vertex triangle $\triangle1$, one of its edge lengths, say
$$\theta^1_{ij}=\cosh^{-1}(2e^{l_{ij}-l_{1i}-l_{1j}}+1),$$
is the remaining edge length of the hyperbolic quadrilateral $Q_{1ij}$ except for $l_{ij}$, $l_{1i}$ and $l_{1j}$. At this point, the three edges $\theta^1_{23}$, $\theta^1_{24}$ and $\theta^1_{34}$ must not satisfy the triangle inequality (note that this triangle is not a face of $\sigma$, so its degeneracy is not considered a two-dimensional degeneracy). In other words, such $l=(l_{12},\ldots,l_{34})$ cannot become the edge length vector of any decorated $1$-$3$ type hyperbolic tetrahedron, resulting in a three-dimensional degeneracy. Such $(\sigma, \{H_2, H_3, H_4\})$ determined by a $l\notin \CL$ is called a virtual $1$-$3$ tetrahedron. It can be viewed as a three-dimensional combinatorial object with two-dimensional geometry and has the same incidence relation and combinatorial structure as the real $1$-$3$ tetrahedra.
Although $\sigma$ corresponds to real geometry on each internal surface, it cannot correspond to any geometric object in $\mathbb{H}^3$, hence the geometric quantities such as $\alpha_{pq}$ and $vol$ have no definition for such $l$.

We summarize the above statement into the following theorem.
\begin{thm}
\label{thm-non-degenerate}
For each $l\in\mathbb{R}^6$, there is no two dimensional degenerations. To be precise, for $\{i, j, k\}=\{2, 3, 4\}$, the two dimensional geometric configurations, i.e. the quadrilateral $p_ip_jv_jv_i$ and the triangle $v_iv_jv_k$ (with corresponding decorations and truncations) always exist, and hence $\theta^1_{ij}$, $\theta^i_{1j}$ and $\theta^i_{jk}$ are always meaningful and can be obtained by the formulas (\ref{Q1ij}), (\ref{cosh-xita}) and (\ref{triangle-XYZ}).
\end{thm}
Moreover, we have the following ``consistency equations".
\begin{prop}
\label{prop-consist}
For all $\{p,q,r\}\subset\{1,2,3,4\}$ and $l\in \mathbb{R}^6$, $\theta^r_{pq}=\theta^r_{pq}(l)>0$. Moreover,
\begin{eqnarray}
\phi_{1i}(l)&=&\frac{\cosh \theta_{ij}^1\cdot\cosh \theta_{ik}^1 - \cosh \theta_{jk}^1}{\sinh \theta_{ij}^1\cdot\sinh\theta_{ik}^1}\\
            &=&\frac{{(\theta_{1j}^{i})}^{2}+{(\theta_{1k}^i)}^{2}-{(\theta_{jk}^i)}^{2}}{2\theta_{1j}^i\cdot\theta_{1k}^i}
\end{eqnarray}
for each $i\in\{2,3,4\}$ and
\begin{eqnarray}
\phi_{jk}(l)&=&\frac{{(\theta_{1j}^{k})}^{2}+{(\theta_{ij}^k)}^{2}-{(\theta_{1i}^k)}^{2}}{2\theta_{1j}^k\cdot\theta_{ij}^k}\\
            &=&\frac{{(\theta_{1k}^{j})}^{2}+{(\theta_{ik}^j)}^{2}-{(\theta_{1i}^j)}^{2}}{2\theta_{1k}^j\cdot\theta_{ik}^j}
\end{eqnarray}
for $\{i,j,k\}=\{2,3,4\}$.
\end{prop}
\begin{proof}
By (\ref{cosh-xita}) and (\ref{triangle-XYZ}), we see $\theta^r_{pq}>0$ for all $\{p,q,r\}\subset\{1,2,3,4\}$.
Substitute these $\{\theta^r_{pq}\}$ into the right-hand side of the above equalities, and further compare with the expressions (\ref{phi-1i}) and (\ref{phi-jk}), we obtain the proof by direct calculations. Since all the calculations are elementary, we omit them here.
\end{proof}

\subsection{Extending the dihedral angles}
\label{extend-angle}
Next we follow the pioneering ideal of Bobenko-Pinkall-Springborn \cite{Bobenko2015}, Luo \cite{Luo2011}, and Luo-Yang in \cite{Luo2018} to extend the definition of the six dihedral angles $\alpha_{pq}(l)$ from $\CL$ to the whole space ${\IR}^6$. Recall the definition of $\phi_{pq}$ in the formulas (\ref{phi-1i}) and (\ref{phi-jk}). For $\{i,j\} \subseteq \{1,2,3,4\}$, let
$$\Omega_{ij}^{\pm} = \{l=(l_{12},\ldots,l_{34})\in {\IR}^6 \mid  \pm \phi_{ij}(l) \geq 1 \},$$
and let
$$X_{ij}^{\pm} = \{(l_{12},\ldots,l_{34})\in {\IR}^6 \mid   \phi_{ij}(l) = \pm 1 \},$$
then $\Omega_{ij}^{\pm}=\Omega_{ji}^{\pm}$, $X_{ij}^{\pm}=X_{ji}^{\pm}$ and $X_{ij}^{\pm}=\partial\Omega_{ij}^{\pm}$. Moreover, by Proposition~\ref{prop 4.4}, we have
\begin{equation}
{\IR}^6\setminus{\CL}=\underset{i\neq j}{\cup} ({\Omega}_{ij}^{+} \cup {\Omega}_{ij}^{-}) ,
\end{equation}

\begin{lemma}\label{lem 4.6}
For $\{ i, j, k, h \}=\{1,2,3,4\}$, we have
\begin{enumerate}
    \item[(1)] ${\Omega}_{ij}^{-}={\Omega}_{kh}^{-}$ and ${\Omega}_{ij}^{+}={\Omega}_{kh}^{+}$.\\[-4pt]
    \item[(2)] $ {\Omega}_{ij}^{-}\cap{\Omega}_{ik}^{-}=\emptyset$, ${\Omega}_{ij}^{+}={\Omega}_{ik}^{-}\sqcup{\Omega}_{ih}^{-}$ and
    ${\Omega}_{ij}^{-}={\Omega}_{ik}^{+}\cap{\Omega}_{ih}^{+}$.\\[-4pt]
    \item[(3)] ${X}_{ij}^{-} = {X}_{kh}^{-}$, ${X}_{ij}^{+}={X}_{kh}^{+}$, ${X}_{ij}^{+}={X}_{ik}^{-}\sqcup{X}_{ih}^{-}$ and ${X}_{ij}^{-}={X}_{ik}^{+}\cap{X}_{ih}^{+}$.\\[-4pt]
    \item[(4)] For each pair $\{i,j\}\subset\{1,2,3,4\}$, we have
    $${\IR}^6\setminus{\CL}={\Omega}_{12}^{-}\sqcup{\Omega}_{13}^{-}\sqcup{\Omega}_{14}^{-}={\Omega}_{ij}^{-} \sqcup {\Omega}_{ij}^{+}.$$
\end{enumerate}
\end{lemma}

\begin{proof}
We first prove the part (1). For all $l\in\mathbb{R}^6$, by the definition of $\phi_{1i}$ in (\ref{phi-1i}), we have
\begin{eqnarray*}
  \phi_{12}(l)&=&
  \frac{2e^{l_{23}-l_{12}-l_{13}}e^{l_{24}-l_{12}-l_{14}}+e^{l_{23}-l_{12}-l_{13}}+e^{l_{24}-l_{12}-l_{14}}-e^{l_{34}-l_{13}-l_{14}}}{2\sqrt{e^{l_{23}-l_{12}-l_{13}}}\cdot \sqrt{1+e^{l_{23}-l_{12}-l_{13}}}\cdot \sqrt{e^{l_{24}-l_{12}-l_{14}}}\cdot \sqrt{1+e^{l_{24}-l_{12}-l_{14}}}}\\
  &=&\frac{2\frac{xy}{d}+x+y-z}{2\sqrt{xy}\sqrt{1+\frac{x}{d}}\sqrt{1+\frac{y}{d}}},
\end{eqnarray*}
where $x=e^{l_{14}+l_{23}}$, $y=e^{l_{13}+l_{24}}$, $z=e^{l_{12}+l_{34}}$ and $d=e^{l_{12}+l_{13}+l_{14}}$.
Hence $l\in\Omega_{12}^{\pm}$ if and only if $\pm\phi_{12}(l)\geq 1$, and if and only if
$$\pm(x+y-z+\frac{2xy}{d})\geq2\sqrt{xy}\sqrt{1+\frac{x}{d}}\sqrt{1+\frac{y}{d}}.$$

On the other hand, by the definition of $\phi_{jk}$ in (\ref{phi-jk}), we have
\begin{eqnarray*}
  \phi_{34}(l)=\frac{x+y-z}{2\sqrt{xy+\frac{xyz}{d}}}.
\end{eqnarray*}
Hence $l\in\Omega_{34}^{\pm}$ if and only if $\pm\phi_{34}(l)\geq 1$, and if and only if
$$\pm(x+y-z)\geq2\sqrt{xy+\frac{xyz}{d}}.$$
By Lemma \ref{lem-elementray-inequa} in the Appendix, the above two inequalities are all equivalent to each other. Hence we obtain ${\Omega}_{12}^{-}={\Omega}_{34}^{-}$ and ${\Omega}_{12}^{+}={\Omega}_{34}^{+}$. Similarly, we have ${\Omega}_{13}^{-}={\Omega}_{24}^{-}$, ${\Omega}_{13}^{+}={\Omega}_{24}^{+}$, ${\Omega}_{14}^{-}={\Omega}_{23}^{-}$ and ${\Omega}_{14}^{+}={\Omega}_{23}^{+}$. By summarizing these properties together and utilizing the symmetry of $\Omega_{ij}$ with respect to subscripts, we finish the proof of part (1).

To prove (2), we just need to show ${\Omega}_{12}^{-}\cap{\Omega}_{13}^{-}=\emptyset$, ${\Omega}_{12}^{+}={\Omega}_{13}^{-}\sqcup{\Omega}_{14}^{-}$
and ${\Omega}_{12}^{-}={\Omega}_{13}^{+}\cap{\Omega}_{14}^{+}$. Because the ideal vertices $v_2$, $v_3$ and $v_4$ are in a similar position, as well as the conclusions in part (1), we obtain the remaining properties in part (2). By Proposition \ref{prop-consist}, we know
$$\phi_{1i}(l)=\frac{\cosh \theta_{ij}^1\cdot\cosh \theta_{ik}^1 - \cosh \theta_{jk}^1}{\sinh \theta_{ij}^1\cdot\sinh\theta_{ik}^1}.$$
It follows that ${\Omega}_{1i}^{-}=\{l\in\mathbb{R}^6:\phi_{1i}(l)\leq -1\}=\{l\in\mathbb{R}^6:\theta_{ij}^1+\theta_{ik}^1\leq\theta_{jk}^1\}$ for $\{i,j,k\}=\{2,3,4\}$.
Hence $l\in {\Omega}_{12}^{-}$ if and only if
\begin{equation}
\label{12-}
\theta_{23}^1+\theta_{24}^1\leq\theta_{34}^1,
\end{equation}
while $l\in {\Omega}_{13}^{-}$ if and only if
\begin{equation}
\label{13-}
\theta_{23}^1+\theta_{34}^1\leq\theta_{24}^1.
\end{equation}
Since all $\theta^r_{pq}>0$, (\ref{12-}) and (\ref{13-}) cannot occur simultaneously, so ${\Omega}_{12}^{-}\cap{\Omega}_{13}^{-}=\emptyset$. Moreover,
\begin{eqnarray*}
{\Omega}_{12}^{+}&=&\{l\in\mathbb{R}^6:\phi_{12}(l)\geq1\}\\
                 &=&\{l\in\mathbb{R}^6:|\theta_{23}^1-\theta_{24}^1|\geq\theta_{34}^1\}\\
                 &=&\{l\in\mathbb{R}^6:\theta_{23}^1\geq\theta_{24}^1+\theta_{34}^1\}\sqcup\{l\in\mathbb{R}^6:\theta_{24}^1\geq\theta_{23}^1+\theta_{34}^1\}\\
                 &=&{\Omega}_{13}^{-}\sqcup{\Omega}_{14}^{-},
\end{eqnarray*}
and
\begin{eqnarray*}
{\Omega}_{12}^{-}&=&\{l\in\mathbb{R}^6: \theta_{34}^1-\theta_{23}^1\geq\theta_{24}^1\}=\{l\in\mathbb{R}^6: \theta_{34}^1-\theta_{24}^1\geq\theta_{23}^1\}\\
                 &=&\{l\in\mathbb{R}^6: |\theta_{34}^1-\theta_{23}^1|\geq\theta_{24}^1\}\cap\{l\in\mathbb{R}^6: |\theta_{34}^1-\theta_{24}^1|\geq\theta_{23}^1\}\\
                 &=&\{l\in\mathbb{R}^6:\phi_{13}(l)\geq1\}\cap\{l\in\mathbb{R}^6:\phi_{14}(l)\geq1\}\\
                 &=&{\Omega}_{13}^{+}\cap{\Omega}_{14}^{+}.
\end{eqnarray*}

Part (3) can be proved similarly. Part (4) is a direct consequence of part (1) and (2). 
\end{proof}

\begin{figure}[htbp]
	\centering
	\subfloat{\includegraphics[width=.450\columnwidth]{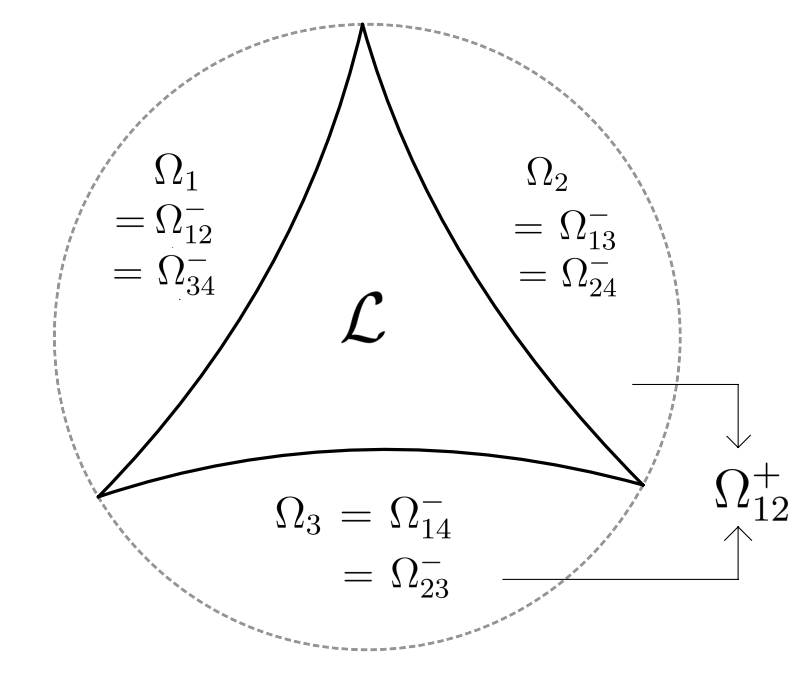}}\hspace{5pt}
	\caption{Decomposition of $\mathbb{R}^6$}
    \label{fig-decom}
\end{figure}

To simplify the notation, we note ${\Omega}_{1}={\Omega}_{12}^{-}$, ${\Omega}_{2}={\Omega}_{13}^{-}$ and ${\Omega}_{3}={\Omega}_{14}^{-}$.
Similarly, we note $ {X}_{1}={X}_{12}^{-}$, ${X}_{2}={X}_{13}^{-}$ and ${X}_{3}={X}_{14}^{-}$.
By Lemma~\ref{lem 4.6}, we know that $X_i = \partial {\Omega}_{i}$ and
 \[
   \partial \CL=\partial({\IR}^6\setminus{\CL})\subseteq \partial {\Omega}_{1} \cup \partial {\Omega}_{2} \cup \partial {\Omega}_{3} = {\sqcup}_{i=1}^3 X_{i} .
 \]
See Figure \ref{fig-decom} for instance.

\begin{prop}\label{prop 4.5}
Let $\partial \CL$ be the boundary of $\CL$ in ${\IR}^6$, then
  $X_{1},\; X_{2},\; X_{3}$ are
  real analytic codimension-1 submanifold of ${\IR}^6$ and
\[
   \partial \CL = {\sqcup}_{i=1}^3 X_{i} .
\]
\end{prop}

\begin{proof}
To simplify the notation, we denote the numberator and the denominator of $\phi_{1i}(l)$ in (\ref{phi-1i}) by $P_i$ and $Q_i$ respectively. For $(l_{12},\ldots,l_{34}) \in {\IR}^6$ and $\{i,j,k\}=\{2,3,4\}$, we have
\[\frac{\partial {\phi}(l_{1i})}{\partial l_{jk}} = -\frac{2}{Q_i} e^{l_{jk}-l_{1j}-l_{1k}} \neq 0 , \]
and then $ \nabla \phi (l_{1i}) \neq 0$. Hence $\phi_{1i}=-1$ is a regular values of $\phi_{1i}$. By the Implicit Function Theorem, $X_i=\phi_{1i}^{-1}(-1)$ is a smooth
co-dimension-1 submanifold of ${\IR}^6$. Since each $\phi_{1i}$ is real analytic in ${\IR}^6$, so the submanifold $X_{i}$ is real analytic.

By Lemma~\ref{lem 4.6}, ${\IR}^6 \setminus  {\CL} = \sqcup_{i=1}^3 {\Omega}_{i}$ is a disjoint union of three six-dimensional submanifolds with boundary. Since $ \partial \CL \subseteq {\sqcup}_{i=1}^3 X_{i} $. We claim that $ X_i \subseteq \partial \CL $. Indeed, for each $l= (l_{12},\ldots,l_{34}) \in X_1 $, we construct a sequence $\{l^{(n)}\} \subseteq \CL $ convergent to $l$ as follows: let $\epsilon_n \to 0^{+}$, define $l^{(n)} = (l_{12}-\epsilon_n,\ldots,l_{34})  \to l $, then for $n$ large enough, by the definition, ${\theta_{23}^1}^{(n)} ,\, {\theta_{24}^1}^{(n)}, \,{\theta_{34}^1}^{(n)} $ satisfy the triangular inequalities, and each $\phi_{1i}(l^{(n)}) \in (-1,1)$, by Proposition~\ref{prop 4.4}, $l^{(n)} \in \CL$ for $n$ large enough.  Therefore, we have $ \partial \CL = {\sqcup}_{i=1}^3 X_{i} $, which completes the proof.
\end{proof}

\begin{figure}[htbp]
	\centering
	\subfloat{\includegraphics[width=.45\columnwidth]{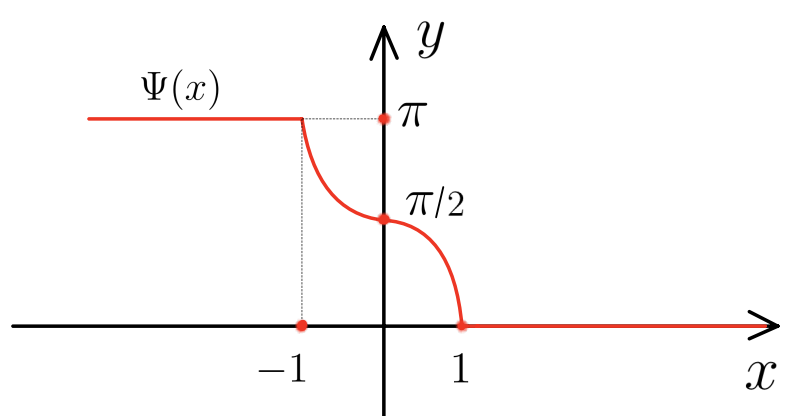}}\hspace{3pt}
	\caption{Function $\Psi(x)$}
    \label{fig-arccos}
\end{figure}

\begin{thm}\label{cor 4.9}
  The dihedral angle function  $\alpha_{ij} : \CL \to \IR$ can be extended continuously to ${\IR}^6$,
  so that its extension, still denoted by $\alpha_{ij}: {\IR}^6 \to \IR $,
  is a constant on each component of $~{{\IR}^6 } \setminus \CL $. To be precise, $(\alpha_{12},\alpha_{13},\alpha_{14},\alpha_{23},\alpha_{24},\alpha_{34})$ equals to
  $(\pi, 0,0,0,0,\pi)$, $(0,\pi,0,0,\pi,0)$ and $(0,0,\pi,\pi,0,0)$ on $\Omega_1$, $\Omega_2$ and $\Omega_3$ respectively.
\end{thm}
\begin{proof}
Extending the definition of those $\alpha_{ij}$ by $\alpha_{ij}\vert_{\Omega_{ij}^{+}}=0$ and
$\alpha_{ij}\vert_{\Omega_{ij}^{-}}=\pi$, for each pair $\{i,j\}\subset\{1,2,3,4\}$, we obtain the proof. Equivalently, one may extend $\alpha_{ij}$ by setting
$$\alpha_{ij}(\cdot)=\Psi(\phi_{ij}(\cdot)),$$
where $\Psi(x)=\arccos(x)$ for $x\in[-1,1]$, $\Psi(x)=\pi$ for $x\leq-1$ and $\Psi(x)=0$ for $x\geq 1$, see Figure \ref{fig-arccos}.
\end{proof}

\subsection{Extending the co-volume function}
\label{section-extend-cov}
The locally convex function co-volume function
$ cov : \CL \to \IR$ satisfies
\[
  \frac{\partial cov}{\partial l_{ij}} = {\alpha}_{ij}
\]
for each pair $\{i,j\}\subset\{1,2,3,4\}$, where ${\alpha}_{ij} : \CL \to \IR $ is the dihedral angle function.
In particular, the differential 1-form $w=\sum\limits_{1\leq i<j\leq 4} {\alpha}_{ij} d l_{ij}$
is closed in $\CL$. Note $\CL$ is simply connected, we can recover the co-volume function $cov$ by $cov(l)=\int^{l}w$.

Since we have extended the dihedral angles $\alpha_{ij}: \CL \to \IR $ to $\alpha_{ij}: {\IR}^6 \rightarrow \IR$,
we define a new continuous 1-form $\mu$ on ${\IR}^6$ by
$$\mu(l) = \sum_{1\leq i<j\leq 4}  \alpha_{ij}(l) \di l_{ij}.$$

The following lemma from Luo~\cite{Luo2011} is very important for latter results.

\begin{lemma}[\cite{Luo2011}]\label{lem 4.11}
  Suppose $U \subseteq \IR^N$ is an open set and $\lambda(l)=\sum_{i} {\alpha}_{i}(x) \di x_{i}$ is
  a continuous 1-form on U,
  \begin{enumerate}
    \item If $ A \subseteq U $ is an open subset bounded by a smooth codimension-1 submanifold of $U$, both $\lambda \vert_A$ and $\lambda \vert_{U\setminus \bar{A}} $ are closed, then $ \lambda $ is closed in $ U $.
    \item If $U$ is simply connected, then $F(x) =\int^x \lambda  $ is a $C^1$-smooth function such that \[ \frac{\partial F}{\partial x_i} = \alpha_i . \]
    \item If $U$ is convex and $A \subseteq U$ is an open subset of $U$ bounded by a codimension-1 real analytic submanifold of $U$ so that $ F \vert_A$ and $F \vert_{U\setminus \bar{A}}$ are locally convex, then $F$ is convex in $U$.
  \end{enumerate}
\end{lemma}

\begin{prop}\label{prop 4.10}
  The differential 1-form $\mu(l)=\sum\limits_{1\leq i<j\leq 4}{\alpha}_{ij}(l) \di l_{ij}$ is closed in $\IR^6$,
  that is, for any Euclidean triangle $\Delta$ in $\IR^6$, $\int_{\partial \Delta} \mu = 0$.
\end{prop}
\begin{proof}
By Theorem \ref{cor 4.9}, each ${\alpha}_{ij}(l) $ is continuous in ${\IR}^6$, hence, the differential 1-form $\mu$ is continuous in ${\IR}^6$. Firstly, it has be showed the restriction $\mu\vert_{\CL} = \sum_{i \neq j} {\alpha}_{ij} \di l_{ij} = \di cov $ is closed. On the other hand, by Proposition~\ref{prop 4.5}, $\CL$ is an open set and bounded by a smooth codimension-1 submanifold in ${\IR}^6$. Let $\overline{\CL}$ be the closure of $\CL$ in ${\IR}^6$, by Theorem \ref{cor 4.9}, $\alpha_{ij}$ is a constant in each connected component of ${\IR}^6 \setminus \bar {\CL} $, hence $\mu$ is closed in ${\IR}^6 \setminus \bar {\CL}$. Using Lemma~\ref{lem 4.11} (a), the differential 1-form $\mu$ is continuous and closed in ${\IR}^6$.
\end{proof}

\begin{cor}\label{cor 4.12}
  The function $cov : \IR^6 \to \IR $ defined by the integral
  \[cov(l) = \int_{(0,\ldots,0)}^{l} \mu + cov(0,\ldots,0) \]
  is a $C^1$-smooth convex function, and satisfies $\nabla_l cov=\alpha$.
\end{cor}

\begin{proof}
By Lemma~\ref{lem 4.11} (3), let $A = \CL$ and $U=\IR^6$, it is easy to show that $cov$ is convex in $\IR^6$. In Proposition~\ref{prop 4.10}, We have proved that $\mu$ is continuous and closed in ${\IR}^6$. By Proposition~\ref{prop 4.10} and Lemma~\ref{lem 4.11} (2), since ${\IR}^6$ is open and simply connected, thus, $cov$ is a $C^1$-smooth function in ${\IR}^6$.

\end{proof}

\section{Global rigidity of $1$-$3$ type hyperbolic tetrahedra}

Let $(M,\CT)$ be a $1$-$3$ type triangulated compact pseudo 3-manifold and let $E=E(\CT)$, $V=V(\CT)$
and $T=T(\CT)$ respectively
be the sets of edges, vertices and tetrahedra in $\CT$. Replacing each 3-simplex in $\CT$ by
a $1$-$3$ type hyperbolic tetrahedron and gluing them along the codimension-1 face by isometries, we obtain a
$1$-$3$ type hyperbolic polyhedral metric on $(M,\CT)$, which is parameterized by the edge length vector $l: E \rightarrow \IR $. This metric is the same as assigning
a sequence of numbers to $E(\CT)$ so that each tetrahedron $\sigma$ becomes a $1$-$3$ type tetrahedron
with assigned numbers as edge lengths.

Denote the space of all $1$-$3$ type polyhedral metrics on $(M,\CT)$ by $\CL(M,\CT)$.
If $l \in \CL(M,\CT)$, its curvature is a map $K_l : E(\CT)\rightarrow \IR $
sending each edge to $2\pi$ minus the sum of dihedral angles
at the edge. The curvature map $K:\CL(M,\CT) \rightarrow {\IR}^{E}$ sends $l$ to $K_l$.

Now, let $w\in \IR^{V(\CT)}$ be a decoration of a $1$-$3$ type tetrahedron. In particular, if $v$ is a hyperideal vertex of some $1$-$3$ type tetrahedron, we define $w(v) = 0$. We define its linearly action on ${\IR}^{E(\CT)}$ as follows:
\[
  (w+x)(vv') = w(v) + w(v') + x(vv')
\]
where $w\in {\IR}^{V(\CT)}$, $x\in {\IR}^{E(\CT)}$, and $vv'$ means the edge with vertices
$v, v'$.

For any two $1$-$3$ type polyhedral metrics $l_1, \;l_2$ on $\CL(M,\CT)$, we can define an equivalent relationship $\sim$
on $\CL(M,\CT)$. We say that $l_1 \sim l_2$ if there exists a change of decoration $ w \in \IR^{V(T)}$ such that
$l_1 = l_2 + w$.
 Furthermore, we define a new space  $\CL(M,\CT) / \sim$ to denote the space $\CL(M,\CT)$ modules
the equivalent relationship $\sim$.

Now, we prove the following rigidity result for decorated $1$-$3$ type hyperbolic polyhedral metrics.

\begin{thm}\label{global rigidity}
Let $M$ be a compact pseudo 3-manifold with a $1$-$3$ type triangulation $\CT$. Then a decorated $1$-$3$ type hyperbolic polyhedral metric on $(M,\CT)$ is determined up to isometry and change of decorations by its curvature.
\end{thm}
\begin{proof}
For each $l\in \CL(M,\CT)/\sim$ and each tetrahedron $ \sigma \in T$, let $l_{\sigma} \in {\IR}^6 $ be the edge length vector of
$\sigma$ in the $1$-$3$ type polyhedral metric $l$. Define the co-volume function
\[
   cov(l)=\sum_{ \sigma \in \CT} cov(l_{\sigma})  .
\]
Since the Hessian matrix of $cov(l_\sigma)$ is positive definite at each point in $\CL(M,\CT) / \sim$, this means the Hessian matrix of $cov$ is
positive definite. Thus, $cov$ is locally strictly convex in $\CL(M,\CT)/ \sim$.

At the same time, the gradient of $cov$ is $2\pi$ minus the curvature map $K_l \in {\IR}^{E}$, that is, $\nabla cov = 2\pi(1,\ldots,1) - K_l $. This could be checked through directly calculation, and Luo first pointed out this when he defined functions $F$ and $H$ in ~\cite{Luo2005}.

By Corollary ~\ref{cor 4.12}, we extend the co-volume function $cov:(\CL(M,\CT)/\sim ) \to {\IR} $ to a $C^1$-smooth convex function, still denoted
by $cov:({\IR}^{E}/\sim) \rightarrow \IR $,
\[
  cov(l) = \sum_{ \sigma \in \CT} cov(l_{\sigma})
\]
where $cov(l_{\sigma})$ is the extended convex function. The convexity of $cov$ follows from the fact that each
summand $cov(l_{\sigma})$ is convex.

Now suppose otherwise that there exist $l_1 \neq l_2 \in \CL({M,\CT})/\sim $ so that
$K_{l_1}=K_{l_2}$. Connecting $l_1$ and $l_2$ in ${\IR}^E$ by the line segment,
 and consider the function $m(t)=cov (tl_1 + (1-t)l_2) $, $ t\in [0,1] $.
Obviously, $ m(t):[0,1]\rightarrow \IR$
is a $C^1$-smooth convex function, and $ m'(t)=\nabla cov \cdot (l_2-l_1) $. Since $\nabla cov(l_i) = (2\pi,\ldots,2\pi) - K_{l_i}$ and $K_{l_1} = K_{l_2}$. Thus $m'(0)=m'(1)$. This means that $m(t)$ must be a linear function in $t$.

On the other hand, since $l_1, l_2 \in \CL({M,\CT})/\sim $, $cov$ should be strictly convex near $l_1$ and $l_2$ along the segment connecting $l_1$ and $l_2$.
Further, $m(t) $ is strictly convex in $t$ near 0 and 1. Thus, there must be $l_1 = l_2$, this contradicts the original assumption.

Therefore, the curvature map $K:\CL(M,\CT)/\sim ~\rightarrow {\IR}^{E}$ must be injective, that is,
for any $l_1 ,\; l_2 \in \CL({M,\CT}) $, if $K_{l_1} = K_{l_2}$, then there exist a change of decoration $w \in \IR^{V(T)}$ such that
$ l_1 = l_2 + w $.
Hence, a $1$-$3$ type hyperbolic polyhedral metric on $(M,\CT)$ is determined up to isometry and change
of decorations by its curvature.

\end{proof}

\section{Volume maximization of angle structure}
\label{section-volume-max}
Let $M$ be a compact pseudo 3-manifold with a $1$-$3$ type closed triangulation $\CT$. Denote $V$, $E$ and $T$ by the sets of vertices, edges and tetrahedra respectively in $\CT$. Let $n$ be the number of tetrahedra in $T$, and we denote the set $T$ by $T=\{\sigma\}$.  Moreover, each $\sigma$ is of $1$-$3$ type, and both the three vertices in the truncated triangle of $\sigma$ are considered as one hyper-ideal vertex. Recall that the sides of the truncated triangles are not considered as the edges in $E$. 
In the following, $e\prec\sigma$ means $e$ is an edge in $\sigma$ and $v\prec e$ means $v$ is a vertex in $e$.


\begin{definition}
An \emph{angle assignment} on $(M,\CT)$ (where $\CT$ is a $1$-$3$ type triangulation) is a function $\alpha:\{(\sigma, e)|\sigma\in T, e\prec\sigma\}\rightarrow \IR_{\geq0}$, so that for each $\sigma\in T$,
\begin{enumerate}
  \item $\sum_{v\prec e}\alpha(\sigma, e)\leq\pi$ for the hyper-ideal vertex $v$ in $\sigma$;
  \item $\sum_{v'\prec e}\alpha(\sigma, e)=\pi$ for each ideal vertex $v'$ in $\sigma$.
\end{enumerate}
If further assume all $\alpha(\sigma, e)>0$ and for each $\sigma\in T$, $\sum_{v\prec e}\alpha(\sigma, e)<\pi$ for the hyper-ideal vertex $v$ in $\sigma$, $\alpha$ is called a \emph{positive angle assignment}. The \emph{cone angle} of an angle assignment $\alpha$ is the function $k_{\alpha}: E\rightarrow{\IR}_{\geq 0}$, $e\mapsto\sum_{e\prec \sigma}\alpha(\sigma, e)$, sending each edge $e$ to the sum of dihedral angles at $e$.
\end{definition}
In the following, $\alpha(\sigma,e)$ is often called the dihedral angle at the edge $e$ in $\sigma$. Further, an angle assignment $\alpha$ may be considered as a point $\alpha\in\IR^{6n}_{\geq0}$, and a positive angle assignment $\alpha$ may be considered as a point $\alpha\in\IR^{6n}_{>0}$. All the angle assignments form a $3n$-dimensional (non-empty) bounded convex polyhedron in $\IR^{6n}_{\geq0}$, the relative interior of which is composed of all the positive angle assignments.

For any given $k\in{\IR}^{E}_{\geq 0}$, which is called a prescribed cone angle, we denote the space of angle assignments
and positive angle assignments with cone angle $k$ by $D_{k}^{*} (M,\CT)$ and $D_{k} (M,\CT)$ respectively. If $\alpha$, $\alpha'$ are two angle assignments,
then $\lambda \alpha+(1-\lambda)\alpha'$ is still an angle assignment for any $\lambda\in[0,1]$. Hence both $D_{k}^{*}(M,\CT)\subset\IR^{6n}_{\geq0}$ and $D_{k}(M,\CT)\subset \IR^{6n}_{>0}$ are convex, and if $D_k(M,\CT) \neq \emptyset $, then $D_{k}^{*}(M,\CT)$ is the compact closure of $D_k(M,\CT)$. For convenience, let's write the vector $(2\pi,\cdots, 2\pi)$ as $2\pi$, then $D_{2\pi}^{*} (M,\CT)$ coincides with the set of all non-negative angle structures on $(M,\CT)$, while $D_{2\pi}(M,\CT)$ coincides with the set of all angle structures $\CA(\CT)$ (see Section \ref{section-introduction}) on $(M,\CT)$.

By Proposition \ref{prop 2.1}, a positive angle assignment $\alpha$ endows each tetrahedron $\sigma\in T$ a $1$-$3$ type hyperbolic tetrahedron so that its dihedral angles are given by $\alpha$. Due to this fact, we can define the volume of a positive angle assignment as the sum of the hyperbolic volume of the $1$-$3$ type tetrahedron $\sigma$, i.e. the volume function $vol:D_{k}(M,\CT)\rightarrow \IR$ is
$$vol(\alpha)=\sum_{\sigma\in T}vol(\sigma).$$

Generally, $D_k(M,\CT)$ is relatively open in some linear subspace in $\IR^{6n}_{>0}$. Moreover, by the work of Leibon \cite{Leibon2002} and Schlenker \cite{Schlenker2002}, the function $vol(\cdot)$ is smooth and strictly concave on $D_k(M,\CT)$.
On the other hand, by the work of Rivin \cite{Rivin2008}, or see Theorem \ref{cor 4.9} and the formula (\ref{volume-alpha}) of volume function,
the function $vol(\cdot)$ can be extended continuously to the compact closure $D_k^{*}(M,\CT)$.

\begin{definition}
Each $l: E \rightarrow {\IR}$ is called a \emph{generalized $1$-$3$ type metric} on $(M,\CT)$.
\end{definition}

For convenience, $l$ is considered as a point in $\IR^{E}$. Recall each $l\in\CL(M,\CT)$ is called a \emph{$1$-$3$ type metric}. If $\sigma\in\CT$ and $e\prec\sigma$, we define the dihedral angle $\alpha(\sigma,e)$ of $l$ to be the corresponding dihedral angles in the generalized hyper-ideal tetrahedron (see Theorem \ref{cor 4.9}) whose decorated edge lengths are given by $l$. Thus we have
\begin{prop}
Each generalized $1$-$3$ type metric $l\in\IR^{E}$ determines an angle assignment $\alpha=\alpha(l)$, while each $1$-$3$ type metric $l\in\CL(M,\CT)$ determines a positive angle assignment $\alpha=\alpha(l)$.
\end{prop}

The main purpose of the rest of the part is to prove the following theorem.


\begin{thm}\label{thm-last-section}
Let $(M,\CT)$ be a closed triangulated pseudo 3-manifold with a $1$-$3$ type triangulation. For any given $k\in{\IR}^{E}_{\geq 0}$ with $D_{k}(M,\CT)\neq \emptyset$, then
  \begin{enumerate}
    \item[(a)] there exists a unique $\alpha\in D_{k}^{*}(M,\CT)$ that achieves the maximum volume.
    \item[(b)] there exists a generalized $1$-$3$ type metric $l\in {\IR}^E$ with cone angle $k$. Moreover, for any generalized metric $l\in {\IR}^E$ with cone angles $k$, the dihedral angles $\alpha(l)$ of $l$ equals to $\alpha$, and hence maximize the volume on $D_{k}^{*}(M,\CT)$.
  \end{enumerate}
\end{thm}
Theorem \ref{thm-last-section} implies Theorem \ref{thm 1.3} with the help of Theorem \ref{cor 4.9}. Using the results in \cite{GJZ1} and \cite{GJZ2}, where the authors showed that every hyperbolic 3-manifold with both cusps and totally geodesic boundaries admits an ideal polyhedral decomposition with each cell either an ideal polyhedron or a partially truncated polyhedron with exactly one truncated face. This suggests that $1$-$3$ type triangulation is likely to occur frequently.
\begin{cor}
Let $M$ be a hyperbolic 3-manifold with both cusps and totally geodesic boundaries, and $\CT$ is a geometric hyper-ideal triangulation of $M$ so that each tetrahedron is either $1$-$3$ type or generalized $1$-$3$ type. If $D_{2\pi}(M,\CT)\neq\emptyset$, then the maximum volume on $D_{2\pi}^{*}(M,\CT)$ is equal to the hyperbolic volume of $M$.
\end{cor}
\subsection{Volume function}
In this section, we return to the study of the geometry in a single $1$-$3$ tetrahedron. Consider a decorated $1$-$3$ type hyperbolic tetrahedron $(\sigma, \{H_2, H_3, H_4\})$ with dihedral angles $\alpha=(\alpha_{12},\ldots,\alpha_{34})\in\CB$,
where $\CB\subset{\IR}^6$ is the space of dihedral angle vectors of all $1$-$3$ type hyperbolic tetrahedra, let $\bar{\CB}$ be its closure, that is
$$\bar{\CB}=\bigl\{(\alpha_{12},\ldots,\alpha_{34})\in{\IR}_{\geq 0}^{6}\,:\,\sum_{j\neq 1}\alpha_{1j}\leq \pi;\;\sum_{j\neq i}\alpha_{ij}=\pi,\,\forall\;i\in\{2,3,4\} \bigr\},$$
where $\alpha_{ij}=\alpha_{ji}$. By (\ref{volume-alpha}), the volume function can be extended continuously to the closed convex polytope $\bar{\CB}$. In the following, each $\alpha\in\bar{\CB}$ is called a generalized dihedral angle (vector) for convenience. We should note that a general $\alpha\in\partial\bar{\CB}$ may not be the extended dihedral angle of any generalized edge length vector $(l_{12}, \cdots, l_{34})\in \IR^6-\CL$. Let
$$\CB_{II}=\{(\pi,0,0,0,0,\pi), (0,\pi,0,0,\pi,0), (0,0,\pi,\pi,0,0)\}$$
and $\CB_{III}=\bar{\CB}-(\CB_{I}\cup\CB_{II})$, where
\begin{equation*}
\CB_I=\bigl\{(\alpha_{12},\ldots,\alpha_{34})\in{\IR}_{\geq 0}^{6}\,:\,\sum_{j\neq 1}\alpha_{1j}<\pi;\;\sum_{j\neq i}\alpha_{ij}=\pi,\,\forall\;i\in\{2,3,4\}\bigr\},
\end{equation*}
then $\bar{\CB}=\CB_I\sqcup\CB_{II}\sqcup\CB_{III}$, $\alpha(\CL)=\CB\subsetneqq\CB_I$ by Proposition \ref{prop 2.1}, and $\alpha(\IR^6-\CL)=\CB_{II}$ by Theorem \ref{cor 4.9}.

By Proposition \ref{prop 2.1}, the corresponding edge lengths $\{l_{ij}\}$ are determined by $\alpha$ up to change of decorations. Those $\{l_{ij}\}$ can be obtained as follows.
Since the dihedral angles are invariant under the change of decorations, 
without loss of generality, we may assume that $l_{12}=l_{13}=l_{14}=0$. Then for $\{i,j,k\}=\{2,3,4\}$,
$$\sinh^2\frac{{\theta}_{ij}^1}{2} = e^{l_{ij} - l_{1i} - l_{1i}} = e^{l_{ij}},$$
it follows $\cosh{\theta}_{ij}^1=1+2e^{l_{ij}}$. On the other hand, since
$$\cosh {\theta}_{ij}^1=\frac{\cos\alpha_{1k}+\cos\alpha_{1i}\cos\alpha_{1j}}{\sin\alpha_{1i}\sin\alpha_{1j}},$$
by direct calculations, we have
\begin{equation}
\label{equ-lij-cos-alpha}
l_{ij}=\ln\frac{1}{2}\Bigl(\frac{\cos\alpha_{1k}+\cos\alpha_{1i}\cos\alpha_{1j}}{\sin\alpha_{1i}\sin\alpha_{1j}}-1\Bigr).
\end{equation}

Let $w \in L(v_1,v_2,v_3)$ be the vector which represent a change of decoration. Recall two metrics $l_1, l_2\in \IR^6$ satisfy the equivalent relationship
$l_1 \sim l_2$ (i.e. differ by a decoration) if there exists a change of decoration
$w\in L(v_1,v_2,v_3)\subset\IR^6$ such that $l_1 = l_2 + w$. Thus the corresponding edge lengths for the given $\alpha=(\alpha_{12},\ldots,\alpha_{34})\in\CB$ is
\[
(\tilde{l}_{12}, \tilde{l}_{13}, \tilde{l}_{14}, \tilde{l}_{23}, \tilde{l}_{24}, \tilde{l}_{34})=(0,0,0,l_{23},l_{24},l_{34})+w.
\]
The diffeomorphism $\tilde{l}:\CB\rightarrow\CL/\sim$, mapping each $(\alpha_{12},\ldots, \alpha_{34})\in\CB$ to the equivalent class
$[(\tilde{l}_{12}, \tilde{l}_{13}, \tilde{l}_{14}, \tilde{l}_{23}, \tilde{l}_{24}, \tilde{l}_{34})]$, may be considered as the inverse map of
$\alpha:\CL\rightarrow\CB$. Each representative element in the equivalent class is still denoted by $\tilde{l}(\alpha)=(\tilde{l}_{12},\ldots,\tilde{l}_{34})\in\CL$, and is called the decorated edge lengths of $\alpha\in \CB$. This will not cause trouble.

Let $S$ be a subset of edges in $\sigma$ satisfying $\{e_{23},e_{24},e_{34}\}\subset S\subset \{e_{ij}\}_{1\leq i<j\leq 4}$, where $e_{ij}$ is an edge in $\sigma$ connecting the two vertices $i$ and $j$, define
\[\CB_S=\{\alpha\in\CB_I\;\vert\;\alpha_{1i} >0\;\text{for}\;e_{1i}\in S\; \text{and}\;\alpha_{1i}=0\;\text{for}\;e_{1i}\notin S\}.\]
\begin{thm}
\label{thm-convex-BS}
The volume function $vol$ is smooth and strictly concave in $\CB_S$.
\end{thm}
\begin{proof}
Recall the definition of the volume function
\begin{eqnarray*}
2vol(\alpha)&=&\Lambda(\alpha_{12})+\Lambda(\alpha_{13})+\Lambda(\alpha_{14})\\
&+&\Lambda \bigl( \frac{\pi-\alpha_{13}-\alpha_{14}+\alpha_{12}}{2} \bigr) +\Lambda \bigl( \frac{\pi-\alpha_{12}-\alpha_{14}+\alpha_{13}}{2} \bigr) \\
&+&\Lambda \bigl( \frac{\pi-\alpha_{13}-\alpha_{12}+\alpha_{14}}{2} \bigr)+\Lambda \bigl( \frac{\pi-\alpha_{12}-\alpha_{13}-\alpha_{14}}{2} \bigr)
\end{eqnarray*}
is extended continuously to $\alpha\in\bar{\CB}$.

If $S=\{e_{23},e_{24},e_{34}\}$, $\CB_S$ is a single-point set, and the conclusion is naturally valid.

If $S=\{e_{ij}\}_{1\leq i<j\leq 4}$, $\CB_S$ is exactly $\CB$, and the conclusion follows from Lemma \ref{hessian-neg-definite}.

If $S$ is composed of five elements, we may assume $S=\{e_{13}, e_{14}, e_{23},e_{24},e_{34}\}$ without loss of generality. Hence $\alpha_{12}=0$, and there are two free variables $\alpha_{13}$ and $\alpha_{14}$ in $\CB_S$ with $\alpha_{13}, \alpha_{14} \in(0,\pi)$ and $\alpha_{13}+\alpha_{14}<\pi$. Using $\Lambda(\theta)+\Lambda(\pi-\theta)=0$, we have
$$2vol(\alpha_{13}, \alpha_{14})=\Lambda(\alpha_{13})+\Lambda(\alpha_{14}) +2 \Lambda \bigl(\frac{\pi-\alpha_{13}-\alpha_{14}}{2} \bigr).$$
By direct calculation, the Hessian matrix $H$ of $2vol$ (with respect to $\alpha_{13}$ and $\alpha_{14}$) is
$$H=-\frac{1}{2}\begin{pmatrix}
           2\cot{\alpha_{13}}+\tan{\frac{\alpha_{13}+\alpha_{14}}{2}} & ~ & \tan{\frac{\alpha_{13}+\alpha_{14}}{2}}\\
           ~&~&~\\
           \tan{\frac{\alpha_{13}+\alpha_{14}}{2}} & ~ & 2\cot{\alpha_{14}}+\tan{\frac{\alpha_{13}+\alpha_{14}}{2}}
\end{pmatrix}.$$
Note $2\cot{\alpha_{13}}+\tan{\frac{\alpha_{13}+\alpha_{14}}{2}}>2\cot{\alpha_{13}}+\tan{\frac{\alpha_{13}}{2}}=\cot\frac{\alpha_{13}}{2}>0$. Similarly,
$2\cot{\alpha_{14}}+\tan{\frac{\alpha_{13}+\alpha_{14}}{2}}>0$. In addition, it's easy to derive
$$\text{det}(H)=\frac{1}{4}\frac{(1+\tan{\frac{\alpha_{13}}{2}}\tan{\frac{\alpha_{14}}{2}})^2}{\tan{\frac{\alpha_{13}}{2}}\tan{\frac{\alpha_{14}}{2}}}>0.$$
Hence $H$ is negative definite, which follows that $vol$ is strictly concave in $\CB_S$.

If $S$ is composed of four elements, we may assume $S=\{e_{14}, e_{23},e_{24},e_{34}\}$ without loss of generality. In this case, there are only one free variable $\alpha_{14}$ in $\CB_S$ with $\alpha_{14}\in(0,\pi)$ and $\alpha_{12}=\alpha_{13}=0$. Then we have
$$2vol(\alpha_{14})=\Lambda(\alpha_{14}) +2 \Lambda \bigl( \frac{\pi-\alpha_{14}}{2} \bigr),$$
which implies
$$\frac{d^2(2vol)}{d\alpha_{14}^2}=-\frac{1}{2}(2\cot{\alpha_{14}}+\tan{\frac{\alpha_{14}}{2}})=-\frac{1}{2}\cot{\frac{\alpha_{14}}{2}}<0$$
by direct calculations. It follows that $vol(\alpha)$ is strictly concave in $\CB_S$.

In conclusion, the volume function is smooth and strictly concave in $\CB_S$.
\end{proof}

\begin{cor}
\label{cor-vol-con-bar-B}
The volume function $vol$ is continuous and strictly concave in $\bar{\CB}$.
\end{cor}

\begin{prop}\label{prop 6.8}
Let $l \in \IR^6$ be the edge length vector of a generalized decorated $1$-$3$ type tetrahedron
with dihedral angle $\alpha(l)\in \CB_{II}$, and let $v \in \IR^6$ so that $\alpha(l) + v \in \CB$, then
\[\lim_{t\to 0^{+}} {\frac{d}{dt}vol\big(\alpha(l)+vt\big)\le -\frac{1}{2}v \cdot l} . \]
\end{prop}
\begin{proof}
We follow the proof strategy of Lemma 6.8 in \cite{Luo2018}. Without loss of generality, we may assume $\alpha=(\pi,0,0,0,0,\pi)$.

Let $f:[0,1] \to \IR$ be the function defined by $f(t)=vol(\alpha(l)+vt)$. By the
concavity of $vol$ and the Mean Value Theorem, we have $f'(t) \le (f(t)-f(0))/t$ for all $t\in (0,1)$.
Since $\alpha(l) \in  \CB_{II}$, $\alpha(l) + v \in \CB$, we know that $f(0)=0$ and $f'(t) < f(t)/t $.

Recall $\CL$ is open in $\IR^6$, $\partial \CL$ is a smooth submanifold, hence for each $l \in \partial \CL$, there exists a
sequence $\{l^{(n)}\}\subseteq \CL$ converging to $l$ and the corresponding dihedral
angles ${\alpha^{(n)}}$ converging to $\alpha(l)$. Since the volume function $vol$ is strictly
concave in $\CB$, we have
$$vol(\beta)-vol(\alpha^{(n)})\le \nabla vol(\alpha^{(n)})\cdot(\beta-\alpha^{(n)})$$
for all $n$ and any $\beta \in \CB$. By the Schl\"{a}fli formula,
$$vol(\beta)-vol(\alpha^{(n)})\le -\frac{l^{(n)}}{2}\cdot(\beta-\alpha^{(n)}).$$
Since $vol$ is continuous on $\bar \CB$ and $vol(\alpha^{(n)}) \to vol(\alpha(l)) = 0$ as $n$ approaches
$\infty$, we know that $vol(\beta)\le -\frac{l}{2}\cdot(\beta-\alpha(l))$.
Now let $\beta=\alpha(l)+vt$, thus $f(t)\le -\frac{t}{2}v \cdot l$.
Combining it with $f'(t) < f(t)/t$, we have $f'(t) \le -\frac{1}{2}v \cdot l$, hence
$\lim_{t\to 0^{+}}f'(t) \le -\frac{1}{2}v \cdot l \;$.

For each $l \notin \bar \CL$, by Proposition \ref{prop 4.4}, there exists an
$m=(l_{12},l_{13},l_{14},l_{23},l_{24}, m_{34})\in \partial \CL$ with $\alpha(m)=\alpha(l)$ and $m_{34} < l_{34}$. By the previous case, we have
$\lim_{t\to 0^{+}}f'(t)\leq-\frac{1}{2}v\cdot m$. Since $\alpha(l) + v \in \CB$, it requires $v_{34}<0$,
so we have $v\cdot l -v\cdot m = v_{34} ( l_{34}-m_{34}) < 0$.
Thus, we also have
\[ \lim_{t\to 0^{+}}f'(t)\le-\frac{1}{2}v\cdot m\leq -\frac{1}{2}v \cdot l .\]

\end{proof}

\begin{prop}\label{prop 6.9}
Assume $\alpha\in \CB_{III}$. If $v \in \IR^6$ satisfies $\alpha+v\in \CB$, then
\[\lim_{t\to 0^+} \frac{d}{dt}vol(\alpha + tv)=+\infty.\]
\end{prop}

\begin{proof}
For any $t\in(0,1)$, $\alpha+tv$ is still in $\CB$, so we can use the Schl\"{a}fli formula, to get
\[\lim_{t\to 0^+}\frac{d}{dt}vol(\alpha+tv)=-\frac{1}{2}\sum_{1\leq i<j\leq 4}v_{ij}\lim_{t\to 0^+}l_{ij}(\alpha+tv).\]

From $\alpha\in B_{III}$, its easy to see $\alpha_{12}+\alpha_{13}+\alpha_{14}=\pi$, $\alpha_{12}=\alpha_{34}$, $\alpha_{13}=\alpha_{24}$, and $\alpha_{14}=\alpha_{23}$.
Since $\alpha+v\in\CB$, then $v_{12}+v_{13}+v_{14}<0$, $v_{23}+v_{24}+v_{34}>0$ and $v\notin L(v_1,v_2,v_3)$. Without loss of generality,
we can assume $l_{12}(\alpha+tv)=l_{13}(\alpha+tv)=l_{14}(\alpha+tv)=0$ for any $t\in(0,1)$.
Now, we only need to show
$$-\frac{1}{2}\Big(\lim_{t\to 0^+}v_{23}\cdot l_{23}(\alpha+tv)+v_{24}\cdot l_{24}(\alpha+tv)+v_{34}\cdot l_{34}(\alpha+tv)\Big)=+\infty.$$
Let's discuss in two cases, first assuming that all the components $\alpha_{ij}$ of $\alpha$ are positive. Since $\alpha+tv\in \CB$, we may use the formula (\ref{equ-lij-cos-alpha}), to get
\begin{eqnarray*}
l_{34}(\alpha+tv)&=&\ln\frac{1}{2}\Big(\frac{\cos(\alpha_{12}+tv_{12})+\cos(\alpha_{13}+tv_{13})\cos(\alpha_{14}+tv_{14})}
{\sin(\alpha_{13}+tv_{13})\sin(\alpha_{14}+tv_{14})} -1 \Big)\\
&=& \ln\frac{1}{2} \Big(\frac{\cos(\alpha_{12}+tv_{12})+ \cos(\pi-\alpha_{12}+t(v_{13}+v_{14}))}{\sin(\alpha_{13}+tv_{13})\sin(\alpha_{14}+tv_{14})}\Big)\\
&=& \ln \frac{1}{2} \Big(\frac{-\sin(\alpha_{12})(v_{12}+v_{13}+v_{14})t+o(t)}{\sin\alpha_{13}\sin\alpha_{14}+O(t)}\Big)\\
&=& \ln t+ C_{34} .
\end{eqnarray*}
Similarly, we have \[ l_{23}(\alpha+tv)=\ln t+ C_{23},\; \; l_{24}(\alpha+tv)=\ln t+ C_{24},\]
thus we obtain
\begin{eqnarray*}
&&-\frac{1}{2}\big(\lim_{t\to 0^+}v_{23}\cdot l_{23}(\alpha+tv)+v_{24}\cdot l_{24}(\alpha+tv)+v_{34}\cdot l_{34}(\alpha+tv)\big)\\&=&
-\frac{1}{2}\lim_{t\to 0^+}(v_{23}+v_{24}+v_{34})\ln t +C =+\infty  .
\end{eqnarray*}
In the second case, some components $\alpha_{ij}$ of $\alpha$ are zero. Since $\alpha\notin \CB_{II}$, there are exactly two $\alpha_{ij}=0$. Without loss of generality, we may assume $\alpha_{12}=\alpha_{34}=0$. In this case $v_{12}>0$. Due to $v_{12}+v_{13}+v_{14}<0$, $v_{12}+v_{23}+v_{24}=0$, $v_{13}+v_{23}+v_{34}=0$ and $v_{14}+v_{24}+v_{34}=0$, we know $2v_{34}=2v_{12}-(v_{12}+v_{13}+v_{14})>0$.
By similar calculation, we may show $l_{34}(\alpha+tv)=2\ln t+ C_{34}$, while both $l_{23}(\alpha+tv)$ and $l_{24}(\alpha+tv)$ are bounded. Hence
$-\frac{1}{2}\lim_{t\to 0^+}v_{34}\cdot \big(l_{34}(\alpha+tv)\big)=+\infty$.
So, in this case we also have
\[\lim_{t\to 0^+} \frac{d}{dt}vol(\alpha+tv) = +\infty .\]
\end{proof}

\subsection{Fenchel dual} 
Let $M$ be a compact pseudo 3-manifold with a $1$-$3$ type closed triangulation $\CT=\{V, E, F, T\}$.
The generalized $1$-$3$ type polyhedral metric on $(M,\CT)$ is parameterized by an edge length vector in ${\IR}^E$. Given $l\in{\IR}^E$ and $\sigma\in T$, denote $l_{\sigma}\in \IR^6$ by the decorated edge length vector of $\sigma$. Define the co-volume function $cov: {\IR}^E \rightarrow \IR$ as
\begin{equation*}
cov(l) = \sum_{  \sigma \in T} cov(l_{\sigma}),
\end{equation*}
where each summand $cov(l_{\sigma})$ is the extended co-volume function introduced in Section \ref{section-extend-cov}.
By Corollary~\ref{cor 4.12},
$cov$ is $C^1$-smooth and convex on $\IR^E$, hence its Fenchel dual
\begin{equation}
  cov^{*}(k) = \sup \{ k\cdot l - cov(l) \; \vert \;l \in  {\IR}^E \}
\end{equation}
is well-defined, convex and lower semi-continuous in ${\IR}^E$. Set
\begin{equation}
  N(\CT) = \{ k \in {\IR}^E \; \vert \; D_{k}^{*} (M,\CT) \neq \emptyset \}.
\end{equation}
For any $k_i\in N(\CT)$ ($i=1,2$), there exists an angle assignments $\alpha_i$ with $k_i=k_{\alpha_i}$. For any $\lambda\in[0,1]$, $\lambda\alpha_1+(1-\lambda)\alpha_2$ is also an angle assignment with cone angle $\lambda k_1+(1-\lambda)k_2$. It follows $\lambda k_1+(1-\lambda)k_2\in N(\CT)$ and hence $N(\CT)$ is a convex subset of $\IR^E$.


\begin{thm}\label{thm 6.10}
Let $ U: {\IR}^E \rightarrow \IR $ be  a function defined by
\begin{equation*}
  U(k) = \left\{
  \begin{array}{cc}
    \min\{-2vol(\alpha) \; \vert \; \alpha\in D_{k}^{*}(M,\CT) \}, & if \; \; k \in N(\CT),\\
    + \infty, & if \; \; k \notin N(\CT) .
  \end{array}
  \right.
\end{equation*}
Then $cov^{*}(k) = U(k)$ for all $k\in {\IR}^E$.
\end{thm}

Before the proof of Theorem \ref{thm 6.10}, we give several propositions.
\begin{prop}
$U$ is convex and lower semi-continuous in $N(\CT)$.
\end{prop}
\begin{proof}
Let $X\in \IR^{6n}_{\geq 0}$ be the space of all angle assignments. $X$ is a $3n$-dimensional non-empty compact convex polyhedron in $\IR^{6n}_{\geq0}$. Set $L:\IR^{6n}\rightarrow \IR^E$ by $L(x)(e)=\sum\limits_{e\prec \sigma}x(\sigma,e)$ which is linear. Note $L(\alpha)=k_{\alpha}$ for each $\alpha\in X$, then $L(X)=N(\CT)$. Set
$$f(\alpha)=-2vol(\alpha)=-2\sum_{\sigma\in T}vol(\alpha_\sigma),\;\alpha \in X.$$
where $\alpha_\sigma=\alpha(\sigma,\cdot)\in \bar{\CB}$ is the dihedral angle vector of $\alpha$ (restricted to $\sigma$). By Corollary \ref{cor-vol-con-bar-B}, $vol(\alpha_{\sigma})$ is continuous and concave on $\bar{\CB}$. Hence $f:X\rightarrow\IR$ is a continuous convex function. Using Lemma \ref{lem-convex-lsc}, we get the conclusion.
\end{proof}

\begin{prop}\label{prop 6.12}
Given $k\in \IR^E$, if $D_k(M,\CT)\neq\emptyset$, then $cov_k:\IR^E \rightarrow \IR$ defined by
  \begin{equation*}
    cov_k (l) = cov(l) - k\cdot l
  \end{equation*}
has a critical point. 
\end{prop}

\begin{proof}
Take any $\alpha \in D_k (M,\CT)$, we rewrite the function $cov_k$ as
  \begin{equation*}
    cov_k (l) = \sum_{  \sigma \in T} (cov(l_{\sigma}) -  \alpha_{\sigma} \cdot l_{\sigma} ) .
  \end{equation*}
For each $\sigma \in T$, we introduce a convex function $ cov_{\sigma,k} : {\IR}^6 \rightarrow {\IR}$, which is defined by
\[cov_{\sigma,k}(x)=cov(x)-\alpha_{\sigma}\cdot x.\]
Let $l(\alpha_{\sigma}) \in \CL$ be an edge length vector of the $1$-$3$ type hyperbolic tetrahedron
whose dihedral angle vector is $\alpha_{\sigma}\in \CB$, 
by Corollary \ref{cor-cov-a*l} and Corollary \ref{cor 4.12}, we know that $l(\alpha_{\sigma})$ is the unique (up to decorations) critical point of the above convex function $cov_{\sigma,k}$.

Note $cov_{\sigma,k}$ is invariant under the linear action of decoration. Consider it as a function defined on the quotient space ${\IR}^6/ L(v_1,v_2,v_3)$, $cov_{\sigma,k}$ is strictly convex near $l_{\sigma}$, so $cov_{\sigma,k}$ is closed in ${\IR}^6/L(v_1,v_2,v_3)$, that is for any $l\in{\IR}^6/L(v_1,v_2,v_3)$,
$\lim_{\vert l \vert \to \infty} cov_{\sigma,k} = +\infty$.

As a consequence, $cov_k(l)$ is convex and closed in ${\IR}^{E}/\sim$.
Hence $cov_k$ has a critical point in ${\IR}^{E}$ and this critical point is unique up do decorations.
\end{proof}

As a consequence of  Proposition~\ref{prop 6.12}, we have

\begin{prop}\label{prop 6.13}
Given $k\in \IR^E$, if $D_k (M,\CT) \neq \emptyset$, then there exists a generalized $1$-$3$ type decorated metric $l\in{\IR}^E$ such that the dihedral angles $\alpha(l)\in D_k^{*}(M,\CT)$.
\end{prop}

\begin{proof}
Each critical point of $cov_k$ provides the desired metric.
\end{proof}

\begin{prop}
The image $ K(\CL(M,\CT)) \cap {(\pi, 2\pi)}^E$ is a convex open polytope in ${\IR}^E$.
\end{prop}
\begin{proof}
The proof is similar with the proof of Proposition 6.14 in Luo-Yang \cite{Luo2018}. Denoting $K(\CT)$ by the convex open (relative to some linear subspace of $\IR^E$) polytope
\begin{equation*}
\{(2\pi,\ldots,2\pi)-k \; \vert \; D_k(M,\CT) \neq \emptyset \}.
\end{equation*}
Next we prove $K(\CL(M,\CT))\cap{(\pi, 2\pi)}^E=K(\CT)\cap{(\pi, 2\pi)}^E$. Note $K(\CL(M,\CT)) \subseteq K(\CT) $, then we have
$K(\CL(M,\CT))\cap{(\pi, 2\pi)}^E\subseteq K(\CT)\cap{(\pi, 2\pi)}^E$. On the other hand, by Proposition~\ref{prop 6.13}, for each $k\in K(\CT)\cap{(\pi, 2\pi)}^E$, there exists a (maybe generalized) $1$-$3$ type metric $l\in {\IR}^E$ such that $K(l)=k$. We prove $l$ is not degenerate, i.e. $l\in \CL(M,\CT)$ as follows. Since $k\in{(\pi, 2\pi)}^E$, the cone angle of $l$ at each edge $e$ is in the range $(0,\pi)$. As a consequence, all the dihedral angles of $\CT$ in $l$ are in the range $(0,\pi)$. So each dihedral angles can not attain $\pi$. Note zero dihedral angle and $\pi$ must occur at the same time (in some $\sigma$), then each dihedral angle can not attain zero. Thus each $\sigma\in T$ is non-degenerate, implying $l\in L(M,\CT)$. Finally we get $ K(\CL(M,\CT)) \cap {(\pi, 2\pi)}^E = K(\CT) \cap {(\pi, 2\pi)}^E$, and the image $ K(\CL(M,\CT)) \cap {(\pi, 2\pi)}^E$ is a convex open polytope in ${\IR}^E$.
\end{proof}

\textbf{Now, we start to prove Theorem \ref{thm 6.10}}.
\begin{proof}
(1) First we show $cov^{*}(k)>C$ for all $C>0$ and all $k\notin N(\CT)$. Consider the compact space $\bar{\CB}^{n}=\bar{\CB}\times\cdots \times\bar{\CB}\subset \IR^{6n}_{\geq 0}$. Since $vol$ is continuous on $\bar{\CB}^{n}$, there exists a constant $C_1>0$ so that $vol(\alpha)\leq C_1$ for all $\alpha\in\bar{\CB}^{n}$. On the other hand, the singleton set $\{k\}$ and $N(\CT)$ are compact and convex sets in ${\IR}^E$, hence for any $C>0$, by the Separation Theorem of Convex Sets,
there exists a $l_0 \in {\IR}^E$, such that
\[
  k\cdot l_0 -c \cdot l_0 > C + 2C_1
\]
for all $ c \in N(\CT)$. In particular, let $c(l_0) \in N(\CT)$ be the cone angle vector of $l_0$, we have
\[
 k\cdot l_0 - c(l_0)\cdot l_0 > C+2C_1 .
\]
Therefore,
\begin{equation*}
  \begin{aligned}
    cov^{*}(k) & \geq k\cdot l_0 - cov(l_0) \\
               & = k\cdot l_0 - c(l_0) \cdot l_0 - 2vol(\alpha(l_0)) \\
               & > C + 2C_1 - 2C_1  =   C ,
  \end{aligned}
\end{equation*}
which means $cov^{*}(k) =\infty$.

(2) Now we prove $cov^{*}(k)=U(k)$ on $N(\CT)$. By Proposition~\ref{prop 6.12}, if $D_k(M,\CT ) \neq \emptyset$, then the function $cov_k$
has a critical point $l$, 
and hence the Fenchel dual of $cov$ is
$$cov^{*}(k)=-cov_k(l)=-2vol(\alpha(l)).$$
For each $\beta\in D_k(M,\CT)$, let $v=\beta-\alpha(l)$, and $v_{\sigma} =\beta_{\sigma} - \alpha_{\sigma}(l)$ for each $\sigma \in T $. Then
$$\sum_{e\prec\sigma}\beta(\sigma,e)=\sum_{e\prec\sigma}\alpha(l)(\sigma,e)=k(e),$$
hence $ \sum_{e\prec\sigma}v(\sigma,e)=0$ for each edge $e$. So we have
$$\sum_{\sigma \in T}v_{\sigma} \cdot l_{\sigma} = \sum_{e\in E} \Big(\sum_{e\prec \sigma}v(\sigma, e)\Big) l(e) = 0.$$
Further, let $T_1$ be the subset of $T$ consisting of (non-degenerate) $1$-$3$ type hyperbolic tetrahedra in $l$ and let $T_2=T\setminus T_1$,
then $-\sum_{\sigma\in T_1}v_{\sigma} \cdot l_{\sigma} = \sum_{\sigma \in T_2}v_{\sigma} \cdot l_{\sigma}$. By the Schl\"{a}fli formula,
\begin{equation*}
\begin{aligned}
\lim_{t\to 0^{+}}\frac{\di}{\di t} vol\big((1-t)\alpha(l)+t\beta\big)&=\lim_{t\to 0^{+}}\frac{\di}{\di t}vol\big(\alpha(l)+tv\big)\\
       &=\sum_{\sigma\in T_2}\lim_{t\to 0^{+}}\frac{\di}{\di t}vol\big(\alpha_{\sigma}(l)+tv_{\sigma}\big)-\frac{1}{2}\sum_{\sigma \in T_1}v_{\sigma} \cdot l_{\sigma}\\
       &=\sum_{\sigma \in T_2}\Big(\lim_{t\to 0^{+}}\frac{\di}{\di t}vol\big(\alpha_{\sigma}(l)+tv_{\sigma}\big)+\frac{1}{2}v_{\sigma}\cdot l_{\sigma}\Big),
  \end{aligned}
\end{equation*}
by Proposition~\ref{prop 6.8}, 
we have
\[
  \lim_{t\to 0^{+}}\frac{\di}{\di t} vol\big(\alpha_{\sigma}(l)+tv_{\sigma}\big) + \frac{1}{2}v_{\sigma} \cdot l_{\sigma}\leq 0,
\]
it follows
\[
  \lim_{t\to 0^{+}} \frac{\di}{\di t} vol\big((1-t)\alpha(l)+t\beta\big)\leq0
\]
for each $\beta\in D_k(M,\CT)$, thus $\alpha(l)$ achieves the maximum volume.
This means $U(k) = -2vol(\alpha(l))$, that is, $cov^{*}$ and $U$ coincide
on the subset $\{k\vert D_k(M,\CT) \neq \emptyset\}$ of $N(\CT)$,
which contains the relative interior of $N(\CT)$.
Since $cov^{*}$ and $U$ are convex and lower semi-continuous, using Lemma \ref{lem-rockafellar}
we have $ cov^{*} (k) = U(k)$ for all $ k \in N(\CT)$.
\end{proof}

\subsection{Proof of Theorem \ref{thm-last-section}}
Now, we prove Theorem \ref{thm-last-section}, which implies Theorem \ref{thm 1.3}.
\begin{proof}
For (a), if it is not true, there exists $\alpha^{(1)}, \alpha^{(2)}\in D_k^{*}(M,\CT)$, both achieve the maximum volume and $\alpha^{(1)}\neq\alpha^{(2)}$.
Connect $\alpha^{(1)}$ and $\alpha^{(2)}$ by  $ L(t) =t \alpha^{(1)} + (1-t)\alpha^{(2)}$, $t\in [0,1]$,
and consider the concave function $f(t) = vol(L(t))$. Because of the maximality of $\alpha^{(i)}$,
the function $f(t)$ must be constant in $[0,1]$.
On the other hand, let  $\alpha_{\sigma}^{(i)} \in \bar{\CB}$ be the restriction of $\alpha^{(i)}$ to $\sigma$ and let
$f_{\sigma}(t) = vol( t \alpha_{\sigma}^{(1)} + (1-t)\alpha_{\sigma}^{(2)} )$ for each $ \sigma \in T$.
Then, by the maximality and Proposition~\ref{prop 6.9}, each $ L_{\sigma}(t) =t  \alpha_{\sigma}^{(1)} + (1-t)a_{\sigma}^{(2)}$ is not in $\CB_{III}$.
So, the interior of the line segment $L_{\sigma}$ lies in $\CB_I\sqcup\CB_{II}$. Since $\CB_{II}$ is discrete in $\bar {\CB}$ and $\alpha^{(1)} \neq \alpha^{(2)}$,
there is at least one ${\sigma}_0 \in T$ such that $L_{{\sigma}_0} \subset \CB_I$.
In fact, the interior of $L_{{\sigma}_0} $ lies in some $\CB_S$ for some subset $S$ of the edges
of $\sigma _0$ (indeed, if $L_{\sigma_0} (t_1) \in \CB_{S_1}$ and $L_{\sigma_0} (t_2) \in \CB_{S_2}$ for some
$t_1, t_2 \in (0,1)$, then the interior of $L_{{\sigma}_0} $ lies in $\CB_{S_1 \cup S_2}$).
By Theorem \ref{thm-convex-BS}, the function $f_{\sigma_0}$ is strictly concave on $\CB_S$, hence
$f=\sum_{ \sigma \in T} f_{\sigma}$ is strictly concave in $(0,1)$. Recall we have shown that $f(t)$ is constant in $[0,1]$. This is a contradiction.  Furthermore, since we have showed that $f_{{\sigma}_0}$ is strictly concave, there is a unique $\alpha \in D_K^{*}(M,\CT)$ that achieves the maximum volume.

For (b), Proposition \ref{prop 6.13} shows that there exists a generalized $1$-$3$ type metric $l\in {\IR}^E$ with cone angle $k$. Further by Proposition \ref{prop 6.12}, $l$ is a critical point of the function $cov_k$. By the following Lemma \ref{lem-cov-vol}, we have $cov_k(l)=2vol(\alpha(l))$.
By Theorem \ref{thm 6.10},
$$cov_k(l)=-cov^{*}(k)=-U(k)=\max\{2vol(\alpha)\, \vert\, a\in D_k^{*}(M,\CT)\}.$$
Hence for each generalized $1$-$3$ type metric $l\in{\IR}^E$ with cone angles $k$, the
dihedral angles $\alpha(l)$ of $l$ is achieved by the maximum volume on $D_k^{*}(M,\CT)$, and $\alpha(l)=\alpha$ by (a).
\end{proof}

\begin{lemma}
\label{lem-cov-vol}
For all $l\in \IR^6$, we have
\begin{equation}
\label{equ-cov-vol-R6}
cov(l)=2vol(\alpha(l))+\sum\limits_{1\leq i<j\leq 4}\alpha_{ij}(l)l_{ij}.
\end{equation}
\end{lemma}
\begin{proof}
By the definition of $cov$, (\ref{equ-cov-vol-R6}) is true on $\CL$. Set $f(l)=cov(l)-\alpha(l)\cdot l$, then on the connected domain $\Omega_1$, $f(l)=cov(l)-\pi(l_{12}+l_{34})$.
Since $\partial cov(l)/\partial l_{ij}=\alpha_{ij}(l)$, we have $\nabla_l f=0$ on $\Omega_1$. Hence $f$ is a constant on $\Omega_1$. For each point $l\in\bar{\Omega}_1\cap\bar{\CL}$, by the continuality, its easy to see $f(l)=2vol(\alpha(l))=2vol(\pi, 0,0,0,0,\pi)=0$. Thus (\ref{equ-cov-vol-R6}) is true on $\Omega_1$. On the domain $\Omega_2$ and $\Omega_3$, we may prove similarly. Thus we get the conclusion.
\end{proof}

At last of this subsection, using the result of volume maximization and angle structures,
we give another proof of the main Theorem~\ref{global rigidity}.

\begin{proof}
Let $l_1,\, l_2 \in \CL(M,\CT)$ be two decorated $1$-$3$ type hyperbolic polyhedral metrics on $(M,\CT)$ with the same cone angles, $k_{l_1} = k_{l_2} = k$.
Then $D_k(M,\CT)\neq \emptyset$ by the assumption. By Theorem \ref{thm-last-section} (a), the dihedral angles
  $\alpha_{l_1}$ and $\alpha_{l_2}$  are both the maximum volume points in
  $D_k^{*}(M,\CT)$. By the uniqueness of the maximum volume point in Theorem \ref{thm-last-section} (a),
   $\alpha_{l_1} = \alpha_{l_2}$. This implies that, for $l_1$ and $l_2$, each pair of corresponding $1$-$3$ type
  hyperbolic tetrahedra are isometric because of Proposition  \ref{prop 2.1}. Therefore, $l_1$ and $l_2$ differs by a change of
  decoration, so we finished the proof.
\end{proof}

\subsection{Proof of Theorem \ref{cor-main-1} and Theorem \ref{cor-main-2}}
We sketch the ideas of proof below.

The part (1) of Theorem \ref{cor-main-1} is known to Lackenby, see Corollary 4.6 in \cite{Lackenby2000}, Theorem 2.1 and Theorem 2.2. in \cite{Lackenby2008}. The basic idea is to consider normal surfaces with combinatorial areas determined by the angle structure $\alpha$ in $(M,\CT)$, the existence of an angle structure $\alpha\in \CA(\CT)$ requires that there is no normal 2-spheres and the only normal tori arise as links of the ideal vertices, hence $M$ is irreducible and atoroidal. By consider properly embedded meridian disc and annulus in normal form, it can be further shown that $M$ is not Seifert fibred. Consider the hyperbolic double $DM$, by Thurston's famous hyperbolization theorem, Lackenby showed that $M$ can be hyperbolized, i.e. $M-\partial_t$ admits a finite-volume hyperbolic structure with totally geodesic boundary.

We remark that the proof of Casson-Rivin's Theorem \ref{thm-casson-rivin} in \cite{Futer2011} can not be generalized to obtain part (2) of our Theorem \ref{cor-main-1}. The main reason is that Thurston's gluing equations (the edge equations and completeness equations)  is not easy to handle in the situation we are considering. To prove the part (2) of Theorem \ref{cor-main-1}, first assume $\alpha\in \CA(\CT)$ is an angle structure, which corresponds to a finite-volume hyperbolic structure with totally geodesic boundary on $M-\partial_t$. Then there is a (non-degenerate) decorated hyperbolic polyhedral metric $l\in \CL$, so that the original ideal triangulation $\CT$ is geometric and $\alpha=\alpha(l)$. By the ``Conversely" part of Theorem \ref{thm 1.3}, $\alpha$ maximizes the volume. Since $\alpha$ is an interior point in $\CA(\CT)$, it is also the critical point of the volume function $vol:\CA(\CT)\rightarrow\IR$. On the contrary, let's assume $\alpha\in \CA(\CT)$ is an angle structure, and $\alpha$ is a critical point of the volume function $vol:\CA(\CT)\rightarrow\IR$. Since $vol(\cdot)$ is smooth and strictly concave on $\CA(\CT)$, and can be extended continuously to the compact closure $\overline{\CA(\CT)}$ (i.e. $D_{2\pi}^{*} (M,\CT)$), it's easy to see $\alpha$ is the maximal point of $vol:\overline{\CA(\CT)}\rightarrow\IR$. Then by Theorem \ref{thm 1.3}, there exists an edge length vector $l\in\IR^E$, so that either case (1) or case (2) in Theorem \ref{thm 1.3} will happen. Since all angles in $\alpha$ are positive, the case (2) will not happen. It follows that the metric $l$ is not degenerate, and hence $l\in\CL$. This $l$ determines a complete finite-area volume hyperbolic structure with totally geodesic boundary on $M-\partial_t$, so that the original triangulation $\CT$ is geometric, and further $\alpha=\alpha(l)$. In other words $\alpha$ corresponds a hyperbolic structure on $M-\partial_t$.

To prove Theorem \ref{cor-main-2}, exactly as in Lackenby's arguments in \cite{Lackenby2008}, one may show any connected closed 2-normal surface in $\CT$ with non-negative Euler characteristic is normally parallel to a toral boundary components of $M$. As a consequence, the proof of Theorem 1.5 in \cite{Garouf}, Theorem 2.5 and Theorem 2.6 in \cite{Kang} can be applied almost word for word to our situation. Hence we omit the details and get Theorem \ref{cor-main-2}.

\section{Appendix}
\subsection{Some elementary lemmas}

\begin{lemma}
\label{lem-elementray-inequa}
For any four positive real number $x$, $y$, $z$, $d$, then
\begin{equation}
x+y+\frac{2xy}{d}+2\sqrt{xy}\sqrt{1+\frac{x}{d}}\sqrt{1+\frac{y}{d}}\leq z
\end{equation}
if and only if
\begin{equation}
x+y+2\sqrt{xy+\frac{xyz}{d}}\leq z.
\end{equation}
Similarly,
\begin{equation}
x+y+\frac{2xy}{d}-2\sqrt{xy}\sqrt{1+\frac{x}{d}}\sqrt{1+\frac{y}{d}}\geq z
\end{equation}
if and only if
\begin{equation}
x+y-2\sqrt{xy+\frac{xyz}{d}}\geq z.
\end{equation}
\end{lemma}
The proof is elementary, we leave it for readers.

\begin{lemma}
\label{lem-convex-lsc}
Suppose $X\subset \IR^n$ is a compact convex set, set $f:X\rightarrow \IR$ is a continuous convex function and $L:\IR^n\rightarrow\IR^m$ is linear. Then $g(y)=\min\{f(x)| x\in X, L(x)=y\}$ is convex and lower semi-continuous on $L(X)$.
\end{lemma}
This is exactly the Lemma 3.6 in Luo-Yang \cite{Luo2018}. The following lemma can be found in Rockafellar's convex analysis textbook (Corollary 7.3.4 \cite{Rockafellar1970}).
\begin{lemma}[\cite{Rockafellar1970}]
\label{lem-rockafellar}
Suppose $X\subset {\IR}^n$ is a closed convex set and $f, g : X \to \IR $ are convex semi-continuous functions. If $f$ and $g$ concide on the relative interior of X, then $f=g$.
\end{lemma}
\subsection{The Schl\"{a}fli formula}
\begin{figure}[htbp]
	\centering
	\subfloat{\includegraphics[width=.7\columnwidth]{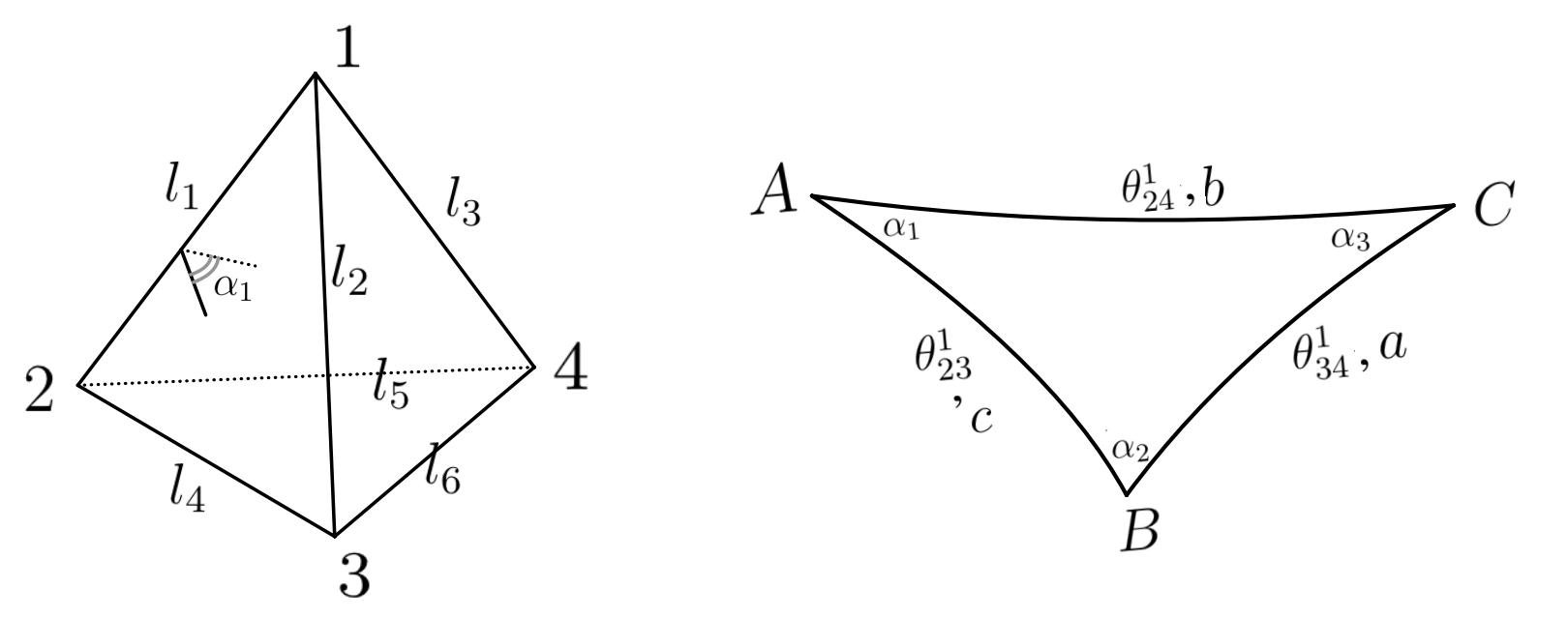}}\hspace{5pt}
	\caption{}
\label{fig-appendix}
\end{figure}

As is shown in Figure \ref{fig-appendix} (left), the edges of a 1-3 type hyperbolic tetrahedron $M_{1|3}$ are relabeled as $l_1$, $l_2$, $\cdots$, $l_6$, and the corresponding dihedral angles are denoted as $\alpha_1$, $\alpha_2$, $\cdots$, $\alpha_6$. As above, the vertex $v_1$ is hyper-ideal and the others are ideal. Consider the volume function $vol(\alpha_1,\cdots, \alpha_6): \CB\rightarrow \mathbb{R}$, we prove the Schl\"{a}fli formula
\begin{equation}
2\di vol+\sum_{i=1}^6l_{i}d\alpha_{i}=0
\end{equation}
directly in this subsection.

Note there are three free variables in these six dihedral angles, such as $\alpha_1$, $\alpha_2$ and $\alpha_3$,
Using the formula (\ref{angle}), we rewrite the volume formula (\ref{volume-alpha}) as the following
\begin{eqnarray*}
2vol(\alpha_{1}, \alpha_{2}, \alpha_{3})&=&\Lambda\Bigl(\frac{\pi-\alpha_{1}-\alpha_{2}-\alpha_{3}}{2}\Bigr)+
\sum_{i=1}^6\Lambda(\alpha_{i}),
\end{eqnarray*}
where $\alpha_4=\frac{\pi+\alpha_{3}-\alpha_{1}-\alpha_{2}}{2}$, $\alpha_5=\frac{\pi+\alpha_{2}-\alpha_{1}-\alpha_{3}}{2}$ and $\alpha_6=\frac{\pi+\alpha_{1}-\alpha_{2}-\alpha_{3}}{2}$.
Note
$$\sum_{i=1}^6l_{i}d\alpha_{i}=\Big(l_1-\frac{l_4+l_5-l_6}{2}\Big)d\alpha_1+\Big(l_2-\frac{l_4+l_6-l_5}{2}\Big)d\alpha_2+\Big(l_3-\frac{l_5+l_6-l_4}{2}\Big)d\alpha_3,$$
we only need to prove
$$2\frac{\partial vol}{\partial \alpha_1}=\frac{l_4+l_5-l_6}{2}-l_1,$$
since the partial derivative formulas for $\alpha_2$ and $\alpha_3$ can also be obtained similarly.
Using
$$\Lambda'(x)=-\ln 2\sin x,$$
its easy to show
$$2\frac{\partial vol}{\partial \alpha_1}=\frac{1}{2}\ln\frac{\cos\frac{\alpha_1+\alpha_2+\alpha_3}{2}\sin\alpha_4\sin\alpha_5}{\cos\frac{\alpha_2+\alpha_3-\alpha_1}{2}\sin^2\alpha_1}.$$
Hence we need to check
$$\frac{\cos\frac{\alpha_1+\alpha_2+\alpha_3}{2}\sin\alpha_4\sin\alpha_5}{\cos\frac{\alpha_2+\alpha_3-\alpha_1}{2}\sin^2\alpha_1}=e^{l_4+l_5-l_6-2l_1}.$$
By the sine law in the Euclidean triangle intersected by the horosphere $H_2$ with the ideal vertex $v_2$, we get
$$\frac{\sin\alpha_4\sin\alpha_5}{\sin^2\alpha_1}=\frac{\theta^2_{14}\theta^2_{13}}{(\theta^2_{34})^2}.$$
Hence we need to check
$$\frac{\cos\frac{\alpha_1+\alpha_2+\alpha_3}{2}}{\cos\frac{\alpha_2+\alpha_3-\alpha_1}{2}}=\frac{\theta^2_{14}\theta^2_{13}}{(\theta^2_{34})^2}\cdot e^{l_4+l_5-l_6-2l_1}.$$
Using $(\theta^2_{34})^2=4e^{l_6-l_4-l_5}$, $\theta^2_{14}=2\sqrt{(e^{l_5-l_1}+e^{l_3})/e^{l_5+l_1}}$, $\theta^2_{13}=2\sqrt{(e^{l_4-l_1}+e^{l_2})/e^{l_4+l_1}}$
and simplifying it, we need to check
$$\frac{1+\cos(\alpha_2+\alpha_3-\alpha_1)}{1+\cos(\alpha_1+\alpha_2+\alpha_3)}=(1+e^{l_1+l_3-l_5})(1+e^{l_1+l_2-l_4}).$$
Further using $\cosh \theta^1_{23}=1+2e^{l_4-l_1-l_2}$ and $\cosh \theta^1_{24}=1+2e^{l_5-l_1-l_3}$, we pneed to check
$$\frac{1+\cos(\alpha_2+\alpha_3-\alpha_1)}{1+\cos(\alpha_1+\alpha_2+\alpha_3)}=
\frac{\cosh\theta^1_{24}+1}{\cosh\theta^1_{24}-1}\cdot\frac{\cosh\theta^1_{23}+1}{\cosh\theta^1_{23}-1}.$$
This is an elementary fact in hyperbolic geometry, i.e. in a hyperbolic triangle $\bigtriangleup ABC$ with edge lengths $a$, $b$ and $c$ which face $A$, $B$ and $C$ respectively, see the right of Figure \ref{fig-appendix}, then one can prove
$$\frac{1+\cos(B+C-A)}{1+\cos(B+C+A)}=\frac{\cosh b+1}{\cosh b-1}\cdot\frac{\cosh c+1}{\cosh c-1}$$
directly by using the standard sine and cosine laws in $\bigtriangleup ABC$.

\begin{rmk}
From the above proof, we see
$$2\frac{\partial vol}{\partial \alpha_1}=-l_1$$
is not correct.
\end{rmk}

\noindent Ke Feng, kefeng@uestc.edu.cn\\[2pt]
\emph{School of Mathematical Sciences, University of Electronic Science and Technology of China, Sichuan 611731, P. R. China}\\[2pt]

\noindent Huabin Ge, hbge@ruc.edu.cn\\[2pt]
\emph{School of Mathematics, Renmin University of China, Beijing 100872, P. R. China}\\[2pt]

\noindent Chunlei Liu, liuchunlei@pku.edu.cn\\[2pt]
\emph{Beijing International Center for Mathematical Research, Peking University, Beijing 100871, P. R. China} \\[2pt]
\end{document}